\documentclass[a4paper,12pt]{amsart}

\usepackage[english]{babel}
\usepackage[latin1]{inputenc}
\usepackage{amsmath}
\usepackage{amsfonts}
\usepackage{amssymb}
\usepackage{enumerate}
\usepackage{geometry}
\usepackage{graphics}
\usepackage[all]{xy}
\usepackage[mathscr]{euscript}
\usepackage{amsthm}
\usepackage{array}
\usepackage{xr}

\newcommand{\pome}{\mathscr{P}}
\newcommand{\cm}{\mathrm{cm}}

\newcommand{\A}{\mathbb{A}}
\newcommand{\B}{\mathbb{B}}
\newcommand{\C}{\mathbb{C}}
\newcommand{\dia}{\diamondsuit}
\newcommand{\diaplus}{\diamondsuit^+}

\newcommand{\uhr}{\!\!\upharpoonright\!\!}
\newcommand{\Talpha}{{T}\uhr{\alpha}}
\newcommand{\Tc}{{T}\uhr{c}}

\DeclareMathOperator{\NF}{succ}

\DeclareMathOperator{\RO}{RO}

\DeclareMathOperator{\mc}{mc}

\DeclareMathOperator{\ot}{ot}

\DeclareMathOperator{\hgt}{ht}

\newtheorem{thm}{Theorem}[section]
\newtheorem{lm}[thm]{Lemma}
\newtheorem{prp}[thm]{Proposition}
\newtheorem{cor}[thm]{Corollary}
\theoremstyle{definition}
\newtheorem{defi}[thm]{Definition}

\newtheorem{question}[thm]{Question}
\theoremstyle{remark}

\begin{document}

% \DOIsuffix{}

% \Volume{}
% \Month{}
% \Year{2010}

% \pagespan{1}{}

% \Receiveddate{}
% \Reviseddate{}
% \Accepteddate{}
% \Dateposted{}

% \keywords{Souslin algebra, Souslin tree, homogeneity}
% \subjclass[msc2000]{03E05, 06E10, 54H05}

\title{Chain homogeneous Souslin algebras}
\author{Gido Scharfenberger-Fabian}
\address{Ernst-Moritz-Arndt-Universit\"{a}t Greifswald\\
Institut f\"{u}r Mathematik und Informatik\\
Walther-Rathenau-Stra\ss e 47\\
D-17487 Greifswald\\
Germany}
\date{\today}

\begin{abstract}
Assuming Jensen's principle $\diaplus$ we construct Souslin algebras all of
whose maximal chains are pairwise isomorphic as total orders, thereby
answering questions of Koppelberg and Todor\c{c}evi\'{c}. %  and 
\end{abstract}
\maketitle

\section*{Introduction}

Souslin's problem was published in the first issue of Fundamenta Mathematicae
(1920, p.223):
\begin{quotation}
Un ensemble ordonn\'{e} (lin\'{e}airement) sans sauts ni lacunes et tel que
tout ensemble de ses intervalles (contenant plus qu'un \'{e}l\'ement)
n'empi\'{e}tant les uns sur les autres est au plus d\'{e}nombrable, est-il
n\'{e}cessairement un continu lin\'{e}aire (ordinaire)?
% \hfill{Probl\`{e}me de M.M.Souslin.}
\end{quotation}
(A (linearly) ordered set without jumps nor gaps and such that each set of its
intervals (containing more than one point) the ones not overlapping onto the
others is at most countable, is it necessarily an (ordinary) linear continuum?)

During the following decades, though still far from being solved, Souslin's
problem was reconsidered by a number of mathematicians who established several
equivalent formulations in terms of trees (\cite{kurepa,miller}), metric
spaces (\cite{kurepa_cr}) or Boolean algebras (\cite{maharam48}).
When in the 1960's the independence of {\sf ZFC} and Souslin's
hypothesis {\sf SH}, which states that the answer to Souslin's
question is ``yes'', was established,
general interest was focussed on the techniques developed for such independence
results, yet there were also some
publications exploring the structural properties of Souslin lines
(\cite{jensen_automorphisms}) and of Souslin trees and algebras
(\cite{jech_automorphisms,jech_simple}). 

Souslin lines can be analyzed through their paritition trees, which are normal
Souslin trees, i.e. normal trees of height $\omega_1$ without uncountable
antichains. There is also an inverse procedure, which generates a Souslin line
from a given Souslin tree. Yet these operations are highly non-canonical,
meaning that one Souslin tree can give rise to up to $2^{\aleph_1}$
non-isomorphic Souslin lines and (in a restricted sense) vice versa.
In the present article we give an example of a Souslin tree to which only one
complete and dense linear order type is associated. 
Moreover, this Souslin tree and its associated Souslin algebra have
interesting homogeneity properties.
And this is were we meet a second, algebraic motivation for the
research presented here.

Homogeneity of a mathematical object generally means that the object locally
looks similar, i.e., that sufficiently small parts of the object are all
structurally equal to each other.
E.g., for a Boolean algebra $B$  homogeneity is defined as the property that
the non-trivial relative algebras $B\uhr a$ for all $a\in B^+$ are isomorphic
to each other.
Of course, many popular structures satisfy stronger homogeneity conditions
than ordinary homogeneity, as e.g., $(\mathbb{R},<)$, the linear order of the
real numbers has for each pair of countable, dense subsets an
auto-homeomorphism that maps one countable, dense subset onto the other.

The variant of homogeneity studied in the present text is \emph{chain
  homogeneity for Boolean algebras}, the property of a Boolean algebra
$B$ that for any two maximal chains $K_0$ and $K_1$ of $B$ (i.e. subsets of
$B$ which
are maximal being linearly ordered by the partial order of $B$) there is an
isomorphism (of linear orders) between $K_0$ and $K_1$.
Note, that such an isomorphism of linear orders does not necessarily extend
to an automorphism of $B$.

In the case of $\sigma$-complete, atomless Boolean algebras it is easy to see
that chain homogeneity implies the c.c.c., so the only possible order types
for maximal chains of such algebras are the real unit interval $[0,1]$ or a
Souslin line with endpoints.
The only known, complete, atomless and chain
homogeneous Boolean algebras are the Cohen-algebras and measure algebras, both
having only separable maximal chains.

The present article gives the affirmative answer to the question of
Sabine Koppelberg whether there are, under appropriate assumptions,
complete and 
atomless chain homogeneous Boolean algebras with inseparable maximal chains.
It is a completely revised version of Chapter 2 of the authors PhD
thesis \cite{diss}. 

The organisation of the text is quite conventional.
In Section \ref{sec:preliminaries} we review the basic definitions concerning
Souslin lines, trees and algebras and fix the notation.
The following three sections introduce the technology we utilize in our
Souslin tree constructions, which consists of relations and
mappings on the tree under construction and results on the extendibility
thereof when adding a new level to the tree.
In Section \ref{sec:hom_for_max} the main results are stated, the proof of
which, a somewhat involved Souslin tree construction, is prepared in Section
\ref{sec:finalprep} and carried out in Section \ref{sec:mainproof}.
The paper closes with a section containing some further remarks and open
problems.

\section{Preliminaries}
\label{sec:preliminaries}
We collect the basic definitions and results and fix the notation
and for later reference.
Most proofs are fairly straght forward and therefore omitted, yet references
are given.
The ambient theory is {\sf ZFC}, Zermelo-Fraenkel set theory with choice;
further hypotheses are explicitely stated.

% Plan: Grundlagen, Wipers, SAE, Kurepa und Anwendungen, Karo und seine
% Anwendung 

\subsection{Boolean algebras}
\label{defi-sect:BA}

Boolean algebras and the most elementary related notions, such as
subalgebras, completeness, atoms will not be defined here,
though we will mention some more or less subtle relationships
between them. The main reference for looking up undefined notions is the
first volume of the \emph{Handbook of Boolean Algebras},
\cite{koppelberg}.

Let $B$ be a Boolean algebra. We call the binary relation on $B$ 
$$x\leq_B y\iff xy=x\iff x-y=0$$
\emph{the natural or the canonical order of} $B$.
A subset $X\subset B$, which is totally ordered by $\leq_B$ is a
\emph{chain} of $B$. A chain $X$ is a \emph{maximal chain} of $B$ if
$X$ is furthermore $\subset$-maximal amongst the chains of $B$.
We call a Boolean algebra all of whose maximal chains are pairwise
order isomorphic \emph{chain homogeneous}.
We let
$$\mc B=\left\{ K\subset B\mid K\text{ is a maximal chain of }B\right\}.$$

An \emph{antichain} of $B$ is a family $X\subset B$ of pairwise
disjoint (i.e. $xy=0$) elements of $B$.
We say that $B$ satisfies (or has) the \emph{countable chain condition}
(or short \emph{c.c.c.}) if every antichain of $B$ is at most countable.

\begin{lm}\label{lm:chainhom-ccc}
If $B$ is $\sigma$-complete and chain homogeneous, then $B$ satisfies the
countable chain condition. 
\end{lm} 
\begin{proof}
Given an uncountable antichain $X$ of $B$ it is easy to construct
well-ordered chains with supremum $1$, $K_0=\{x_0<_B x_1<_B\ldots\}$ of order type $\omega$ and
$K_1=\{y_0 <_B y_1 <_B\ldots\}$ of order type $\omega_1$. Now for any pair $K,\,K'\in\mc B$ with
$K_0\subset K$ and $K_1\subset K'$ there should be an isomorphism $\varphi:K\to K'$
but then $n\mapsto\min\{\alpha\mid \varphi(x_n)\leq_B y_\alpha\}$
would give a cofinal countable sequence in $\omega_1$.
\end{proof}
% Now let $B$ be a complete Boolean algebra.
% By saying that $A$ is a \emph{complete subalgebra of} $B$ we mean that
% $A$ is a subalgebra of $B$ which is a complete Boolean algebra,
% and it computes the same infinite sums and products as $B$ does.
% I.e., for all $M\subseteq A$ we have
% $${\sum}^A M={\sum}^B M$$
% and the analogous equation for products holds as well.
For any subset $X$ of a complete Boolean algebra $B$ we call
$$\langle X\rangle^\cm:=
\bigcap\left\{A\mid X\subseteq A\text{ and }A\text{ is a complete subalgebra of }B\right\}$$
the \emph{subalgebra of $B$ that is completely generated by $X$} and maybe write
$\langle X\rangle^\cm_B$ if $B$ is not clear from the context.
If the superscript ``$\cm$'' is omitted, i.e., by $\langle X\rangle_B$,
we denote the intersection of all subalgebras containing $X$ as a
subset (and not only the complete subalgebras) and call $\langle X\rangle$ the subalgebra of $B$, that
is \emph{finitarily} generated by $X$.

Note that given an arbitrary Boolean algebra $A$ which is a subalgebra of a complete
Boolean algebra $B$, the \emph{(Dedekind) completion} $\overline{A}$%=\RO A^+$
(cf. \cite[Section 4.3]{koppelberg}),
which is the unique complete Boolean algebra
containing $A$ as a dense subalgebra, 
is not necessarily isomorphic to $\langle A\rangle^\cm_B$.
A \emph{regular subalgebra}\label{reg_subalg}
$A$ of $B$ is a subalgebra,
such that for all $M\subset A$ we have
${\sum}^B M={\sum}^A M$ if the latter sum exists.
(This definition also extends to \emph{regular sublattices}.)
% Here $A$ need not be complete and
% only the sums that exist in $A$ are considered.
A well-known fact is, that $A$ is a regular subalgebra of $B$ if and
only if $A$ is a dense subset of $\langle A\rangle_B^\cm$.
In this case $\langle A\rangle_B^\cm$ is indeed isomorphic to $\overline{A}$,
%:=\RO(A^+,<_B)$,
the Dedekind completion of $A$ (cf. \cite[Prop.4]{cohen_algebras}).

A Boolean algebra $B$ is $\aleph_0$\emph{-distributive} if for every family
$(a_{ij})_{i\in\omega,j\in J}$ with an index set $J$ of arbitrary size,
the following equation holds:
$$\prod_{i\in\omega}\sum_{j\in J} a_{ij}=
\sum\left\{\prod_{i\in\omega}a_{if(i)}\mid f\in{}{^\omega J}\right\}.$$
We do not need more specific concepts of distributivity.

% An antichain is a \emph{maximal antichain} or a \emph{partition (of unity)}
% if it is $\subset$-maximal,
% and if $X,Y$ are maximal antichains we say that $X$ \emph{refines} $Y$
% if for all $x\in X$ there is a $y\in Y$ above $x$, i.e. $x\leq_B y$.
% With this notion at hand we can state
% a useful characterisation of of $\aleph_0$-distributivity:
% A Boolean algebra $B$ is $\aleph_0$-distributive if and only if for every countable
% family $(Y_n)$ of maximal antichains of $B$ there is a common
% refinement $X$, i.e. a maximal antichain that refines all the $Y_n$
% (cf. \cite[Proposition 14.9]{koppelberg}).

We frequently consider the \emph{regular open algebra} $\RO X$ of some
topological space $X$.
A subset $U$ of $X$ is regular open
if the interior of the closure of $U$ is equal to $U$ and
$\RO X$ is the set of all regular open subsets of $X$.
The regular open algebra of any space $X$ is a complete
Boolean algebra
but it is in general not a subalgebra of $\pome(X)$
from which the operations are modified by taking the
regularisations (cf.~\cite[Theorem 1.37]{koppelberg}).
If the topology of $X$ is not specified,
then one of following applies:
\begin{itemize}
\item $X$ is a linear order; then $X$ carries the order topology
(as in the following proposition),
\item $X=T$ is a tree; then $T$ carries the partial order topology
which is generated by all subsets of the form $T(t):=\{s\in T\mid t\le_T s\}$
for some $t\in T$ (a declaration of our tree notation follows below).
\end{itemize}

\subsection{Normal trees}

% \subsection{Normal trees}
% \label{sec:normaltrees}

A \emph{tree} is a well-founded partial order $(T,<_T)$ where
the set of predecessors $\{s\mid s<_T t\}$
is totally (and thus well-) ordered by $<_T$ for all $t\in T$.
The elements of a tree are called \emph{nodes}.
The \emph{height} of the node $t$ in $T$ is the order
type of the set of its predecessors under the ordering of $T$,
$\hgt_T(t):=\ot(\{s\mid s<_T t\},<_T)$.
For an ordinal $\alpha$ we let $T_\alpha$ denote
the set of nodes of $T$ with height $\alpha$.
If $\hgt_T(s)>\alpha$ we let $s\uhr\alpha$ be the unique
predecessor of $s$ in level $T_\alpha$.

The \emph{height of a tree} $T$, $\hgt T$,
is the minimal ordinal $\alpha$ such that $T_\alpha$ is empty.
An \emph{antichain} is a set of pairwise incomparable nodes of $T$,
so for $\alpha<\hgt T$,
the level $T_\alpha$ is an antichain of $T$.

If $c$ is a subset of the height of a tree $T$, let $\Tc$ denote the tree
that consists of all the nodes of $T$ whose height lies in $c$
together with the inherited tree order $<_T$:
$$T\uhr c=\bigcup_{\alpha\in c}T_\alpha,\quad s<_{T\upharpoonright c}t\iff s<_T t.$$
% In this context it is worth
% noting that $\hgt_{T\upharpoonright c}(t)=\ot(c\cap\hgt_T(t))$ for all $t\in\Tc$
% which is possibly smaller than $\hgt_T(t)$.

For a node $t\in T$ we let $\NF(t)$ be the set of $t$'s immediate successors.
By $T(t)$ we denote the set of all successors of $t$ in $T$ (including $t$)
which forms a tree with root $t$ under the ordering inherited from $T$.
Nodes, that do not have $<_T$-successors, are called \emph{leaves}, and
$T$ is called $\kappa$-branching, $\kappa$ a cardinal, if all nodes of $T$
have exactly $\kappa$ immediate successors, except for the leaves.

A \emph{branch} is a subset $b$ of $T$ that is linearly ordered by $<_T$ and
closed downwards, i.e. if $s<_T t\in b$ then $s\in b$.
For $\alpha<\ot(b,<_T)$ we let $b\uhr\alpha$ be the unique element of $T_\alpha\cap b$ and so 
extend the similar notation for nodes and their predecessors to branches and their elements.
(For a normal tree this is just natural,
since the nodes can be identified with the branches leading to them.)
A branch is \emph{maximal} in $T$
if it is not properly contained in any other branch of $T$. 

Under the notion of a \emph{normal} tree we subsume the following
four conditions:
\begin{enumerate}[(i)]
\item there is a single minimal node called the \emph{root};
\item each node (except for the top level nodes) has at least two immediate
  successors; 
\item each node has successors in every higher non-empty level;
\item branches of limit length have unique limits (if they are extended in the tree),
i.e., if $s,t$ are nodes of $T$ of limit height
whose sets of predecessors coincide, then $s=t$.
\end{enumerate}

Whenever we consider a mapping $\varphi:T\to S$
between trees and call this mapping a \emph{tree homomorphism} we mean that
$\varphi$ carries $<_T$ to $<_S$ and respects the height function:
$\hgt_S(\varphi(s))=\hgt_T(s)$ for all $s\in T$. 

A tree $T$ is said to be \emph{homogeneous},
if for all pairs $s,t\in T$ of the same height there is
a tree isomorphism between $T(s)$ and $T(t)$,
the trees of nodes in $T$ above $s$ and $t$ respectively.

\subsection{Souslin lines, trees and algebras}

% \subsubsection{Souslin lines}
It is well-known that the existence statements of Souslin algebras, lines
and trees are mutually equivalent, yet independent of {\sf ZFC}: they all
exist if and only if Souslin's hypothesis {\sf SH} fails.

A \emph{Souslin line} is a total order directly witnessing the failure of  {\sf SH},
ie., it is a complete, dense linear order that satisfies the 
countable chain condition but is not separable,
i.e., it has no countable dense subset, but all families
of pairwise disjoint open intervals are countable.

In this definition the possibility of a Souslin line
having a separable non-trivial interval is included.
We will, however,
never consider a Souslin line with a non-trivial separable interval.

It is a standard result, that
every Souslin line $L$ has a dense subset of cardinality $\aleph_1$
and therefore the cardinality $|L|=2^{\aleph_0}$.

% \subsubsection{Souslin trees}

In general,
a \emph{Souslin tree} is a tree $T$ of height $\omega_1$
such that
every family of pairwise incomparable nodes
and every branch of $T$ are at most countable.
We will only consider normal Souslin trees, where absence of uncountable
antichains already implies that the tree has no cofinal branch.

Souslin lines and Souslin trees are tightly related by the following
operations (also cf. \cite{todorcevic_trees_and_orders} or
\cite[9.14]{jechneu}). 
If we are given a Souslin tree $T$ and a total order of its supporting set
(or just total orderings of the sets of the form $\NF(t)$ for all $t\in T$)
we can equip the set of all maximal branches of $T$
with the lexicographical ordering.
The resulting linear order is dense, c.c.c.\ and inseparable, i.e., its
Dedekind completion is a Souslin line.

Now let $L$ be a given Souslin line (without separable intervals).
Then every normal partition tree of $L$ is a Souslin tree.
A partition tree of a linear order $L$ is a tree whose nodes are open
intervals of $L$ such that the union of each level is dense in $L$ and whose
tree ordering is given by the inverted inclusion.

Anyway, in general these operations are highly non-canonical; one Souslin
tree may be associated to many non-isomorphic Souslin lines and vice versa.
To unify considerations one uses the notion of Souslin algebras.
% \subsubsection{Souslin algebras}

A \emph{Souslin algebra} is a
complete, atomless, c.c.c.\ and $\aleph_0$-distributive
Boolean algebra\footnote{In other contexts than ours, $\kappa$-Souslin algebras are defined as
complete, atomless, $\kappa$-c.c. and $(<\kappa)$-distributive
Boolean algebras and $\kappa$ can be any
uncountable cardinal.
In this notation the objects of our consideration are called
$\aleph_1$-Souslin algebras.
But by Lemma (\ref{lm:chainhom-ccc}) above,
these higher Souslin algebras fo $\kappa>\aleph_1$
always have maximal chains of distinct order types
and are therefore never chain homogeneous.}.
We will use the following non-standard denotations:
A Souslin algebra that has a set of complete generators of size $\aleph_1$ will
be called \emph{small}; otherwise it is \emph{big}.
Small Souslin algebras always have $2^{\aleph_0}$ elements, while a famous
result of Solovay states that the cardinality of (big) Souslin algebras is 
at most $2^{\aleph_1}$ (cf.\~{\cite[Thm.30.20]{jechneu}}).

We use letters $\A,\B,\mathbb{C}$ to denote Souslin algebras.
In a Souslin algebra $\B$, a complete subalgebra $\A$ is atomic if and only if
$\A$ is completely generated by some countable subset $X$ of $\B$
(cf. \cite[Prop.14.8]{koppelberg}).
In the other extreme,
since distributivity and the c.c.c.
are handed down to complete subalgebras,
the atomless, complete subalgebras of a
Souslin algebra $\B$ are Souslin algebras as well.

We now review the connections between Souslin algebras and Souslin lines and
trees respectively, beginning with the lines.

\begin{lm}[{\cite[\S 14, Exrc.5]{koppelberg}}, for a proof cf. {\cite[2.1.2-5]{diss}}]
\begin{enumerate}[a)]
\item The maximal chains of a Souslin algebra
$\B$ are Souslin lines with endpoints and without
non-trivial separable intervals.
\item Let $L$ be a Souslin line without separable intervals.
The regular open algebra of $L$ is a Souslin algebra.
Furthermore,
if $L$ has endpoints then $\RO L$ has a maximal
chain $K$ which is isomorphic to $L$
such that $\langle K\rangle^\cm=\RO L$.
\item Let $\B$ be a Souslin algebra and $K\in\mc\B$.
Then $\langle K\rangle$, the subalgebra of $\B$ that is finitarily generated by
$K$, is a regular subalgebra of $\B$.
In particular we have
$$\langle K\rangle^\cm\cong\overline{\langle K\rangle}\cong\RO K.$$
\end{enumerate}
\end{lm}

It is well known that for a normal Souslin tree $T$, its regular open algebra
$\RO T$ is a Souslin algebra and $T$ can be densely embedded in $\RO T$.
On the other hand, in every Souslin algebra $\B$ that has
a family of complete generators of cardinality $\aleph_1$, i.e., $\B$
is a \emph{small} Souslin algebra, 
there is a dense subset $T$ of $\B\setminus\{0\}$
such that $(T,>_\B)$ is a normal Souslin tree
(note that $<_T$ is $>_\B$) that is regularly embedded in $\B$ and
whose regular open algebra is therefore
isomorphic to $\B$ (see e.g. \cite[Thm.14.20]{koppelberg}).
We will use the following convention established in \cite{devlin-johnsbraten}.
Let $T$ be a subset of a Souslin algebra $\B$.
Then $T$ is said to \emph{souslinise} $\B$ or
to be a \emph{Souslinisation} of $\B$ if $T$ is dense in $\B$ and
becomes a Souslin tree under the reverse Boolean order of $\B$. 
Every level $T_\alpha$ of a Souslinisation $T$ of $\B$ is a partition of
unity in $\B$ and
taking limits in $T$ is simply the evaluation of the
corresponding infinite product in $\B$:
If $t\in T_\alpha$ is a limit node,
let $t_\gamma$ for $\gamma<\alpha$ be
the unique $>_\B$-predecessor of $t$ of height $\gamma$.
Then $t=\prod\{t_\gamma\mid\gamma<\alpha\}$.
The Souslinisation is unique up to the elimination of a non-stationary set of
levels: Given an isomorphism $\varphi$ between two  Souslin algebras with
Souslinisations $T_1$ and $T_2$ respectively, there is a club
$C\subseteq\omega_1$, such that $\varphi$ becomes an isomorphism of trees when
restrictied to $T_1\uhr C$ and $T_2\uhr C$.
For a proof of this well-known result, cf.~\cite[Lemma
VIII.9]{devlin-johnsbraten} or \cite[Lemma 25.6]{jech}.
% \begin{lm}[Restriction Lemma]
% \label{lm:souslinisation}
% Let $\B_0,\,\B_1$ be Souslin algebras
% with Souslinisations $T_0$ and $T_1$ res\-pectively.
% Then a maping $\varphi:\B_0\to\B_1$ is an isomorphism just in case that
% there is a closed unbounded subset $C$ of $\omega_1$ such that
% the restriction of $\varphi$ to $T_0{}\uhr{}C$ is a tree isomorphism onto $T_1{}\uhr C$.
% \end{lm}

% \begin{rem}
% It is easily seen that each tree automorphism of the Souslinisation $T$
% extends to a Boolean automorphism on $\B$,
% yet not every automorphism of $\B$ needs to restrict to $T$.
% Anyway, the last theorem implies that there must be some $c$ club in $\omega_1$
% such that the Boolean automorphism restricts to $T\uhr c$.
% Therefore, we call Lemma \ref{lm:souslinisation} \emph{the Restriction Lemma
% (for isomorphisms between Souslin algebras)}.
% \end{rem}

A Souslinisation $T$ of $\B$ provides
a natural stratification of $\B$ by
countably generated, complete and therefore atomic subalgebras.
Fix a Souslinisation $T$ of $\B$ and let for $\alpha<\omega_1$
$$\B^\alpha:=\langle T_\alpha\rangle^\cm.$$
Note that for all $\alpha<\omega_1$ we have $\B^\alpha\cong\pome(\omega)$
and $T\uhr(\alpha+1)\subset\B^\alpha$.
Clearly the sequence of the $\B^\alpha$ is increasing.
To show $\B=\bigcup_{\alpha<\omega_1}\B^\alpha$,
pick $a\in\B\setminus\{0\}$.
There is a maximal pairwise disjoint subset $A$ of $T$ of elements $\leq_\B a$,
so $\sum A=a$.
By the countable chain condition $A$ must be countable and therefore a subset of
$\B^\alpha$ for some $\alpha<\omega_1$.
We finally note, that the sequence of the $\B^\alpha$ is not continuous.
For a countable limit ordinal $\alpha$ we have
$$\bigcup_{\gamma<\alpha}\B^\gamma\lneq
\left\langle\bigcup_{\gamma<\alpha}\B^\gamma\right\rangle^\cm=\B^\alpha.$$

\subsection{Jensen's principle $\diaplus$}

As for many Souslin tree constructions in the literature,
we will assume diamond-principles.

\begin{defi}
\begin{enumerate}[a)]
\item The sequence $(R_\alpha)_{\alpha<\omega_1}$ of sets
$R_\alpha\subseteq\alpha$ for $\alpha<\omega_1$
is a $\dia$-sequence if
for all $X\subseteq\omega_1$ there is some stationary set $s\subseteq\omega_1$
such that, for all $\alpha\in s$ we have $X\cap\alpha=R_\alpha$.
The statement ``There is a $\dia$-sequence.''
will be denoted by $\dia$.
\item The sequence $(S_\alpha)_{\alpha<\omega_1}$ of countable sets
$S_\alpha\subseteq\pome(\alpha)$ is a $\diaplus$-sequence if
for all $X\subseteq\omega_1$ there is some club set $c\subseteq\omega_1$
s.t., for all $\alpha\in c$ we have $X\cap\alpha,\, c\cap\alpha\in S_\alpha$.
The statement
``There is a $\diaplus$-sequence.''
will be denoted by $\diaplus$.
\end{enumerate}
\end{defi}

% If we delete the requirement that
% $c\cap\alpha$ be in $S_\alpha$ for $\alpha\in c$
% in the definition of the $\diaplus$-sequence,
% then we define a $\diastar$-sequence.
% This will not be used here.
Since $\diaplus\rightarrow\dia$
we only need to assume $\diaplus$ in our
statements even though for the sake of convenience we will use both,
a $\diaplus$-sequence and a $\dia$-sequence, at the same time.

\section{Generating chains and wipers}

For the representation of a given maximal chain $K$ of a Souslin algebra $\B$
we will use the approximations
$K_\alpha=K\cap\B^\alpha$ with respect to some fixed Souslinisation
$T$ which are chains in $\B^\alpha$, yet not necessarily maximal ones.
In the simple case of power set algebras or their isomorphic copies,
such as the $\B^\alpha$,
we can give a nice characterisation of maximal chains.

\begin{prp}[{folklore, \cite[2.1.10]{diss}}]
\label{prp:max_chn_lin}
Let $X$ be a set.
Given a $\subset$-chain $K$ of $\pome(X)$,
define the quasi-ordering $\leq_K$ on $X$ by
$$x\leq_K y\leftrightarrow(\forall u\in K)(y\in u\Rightarrow x\in u)$$
and let $\overline{K}:=\{\sum M,\,\prod M\mid M\subseteq K\}$ be the
completion of the linear order $(K,\subset)$ in $\pome(X)$.
Then the following statements are equivalent:
\begin{enumerate}[(i)]
\item $\overline{K}$ is a maximal chain of $\pome(X)$;
\item $\leq_K$ is a linear ordering of $X$;
\item $\langle K\rangle^\cm_{\pome(X)}=\pome(X)$.
\end{enumerate}
In particular, there is a bijective association between the total orders on
$X$ and the maximal chains of $\pome(X)$.
\end{prp}
% \end{rem}

\begin{defi}
A maximal chain $K$ of a Souslin algebra $\B$ that completely
generates $\B$ will be called a \emph{generating chain}.
\end{defi}
We stress that a chain $K$ that satisfies $\langle K\rangle^{\cm}=\B$
will not be called generating unless it is maximal.
\begin{lm}[cf. {\cite[2.2.2]{diss}}]
\label{lm:genchaintree}
Given a generating chain $K\subseteq\B$ the set
$$C:=\{\alpha<\omega_1:K\cap \B^\alpha\in\mc(\B^\alpha)\}=
\{\alpha\in\omega_1:\langle K\cap\B^\alpha\rangle^\cm=\B^\alpha\}$$
is closed and unbounded in $\omega_1$.
\end{lm}
The proof of the last lemma works by a straight forward catch-up
argument using the c.c.c. and the distributive law of $\B$.

Combining Proposition \ref{prp:max_chn_lin} and Lemma
\ref{lm:genchaintree}, we find a means for the representation of
generating chains of a Souslin algebra as a relation on one of its
Souslinizations that we now fix in the notion of
``wiper''. Requirement (ii) below is of technical importance when it
comes to the construction of isomorphisms between wipers.

\begin{defi}\label{defi:wipers}
\begin{enumerate}[a)]
\item Let $T$ be a normal tree with countable levels
of arbitrary height $\alpha\leq\omega_1$
and let $C$ be a subset of $\alpha$.
A \emph{wiper of total orders} (or
more convenient a \emph{wiper}) on $T{}\uhr C$ is a family
$W=\langle <_\gamma\mid\gamma\in C\rangle$
of total orders $<_\gamma$ on $T_\gamma$
such that
\begin{enumerate}[(i)]
\item $W$ respects the tree order of $T$: for all $\beta,\gamma\in C$,
$\beta<\gamma$ and $s,t\in T_\beta$, $s',t'\in T_\gamma$ we have
$$s<_T s' \wedge t<_T t' \wedge s<_\beta t\Rightarrow s'<_\gamma t',$$
\item if $\beta,\gamma\in C$, $\beta<\gamma$, and $s\in T_\beta$,
then the set $I_{s,\gamma}=\{t\in T_\gamma\mid s<_T t\}$ of successors of $s$ in
level $\gamma$ is ordered densely without endpoints by $<_\gamma$.
\end{enumerate}
\item If $T$ of height $\omega_1$ is a Souslinisation of $\B$ and $W$ is a
  wiper on $T$, we say that
$$K_W:=\big\{\sum M\mid(\exists\alpha\in \omega_1)M\subseteq T_\alpha,\,
M\text{ is an initial segment of }<_\alpha\big\}$$
is the subset of $\B$ which is \emph{induced} by $W$.
\end{enumerate}
\end{defi}

% Only the first condition in this definition is inspired by Lemma \ref{lm:genchaintree}
% while the second is a useful standardisation requirement
% whose technical importance will become clear in the proof of
% Proposition \ref{prp:extendwiperisos} where we aim at building isomorphisms between wipers.
To illustrate the definition and justify its definiendum
imagine the Souslin tree $T\uhr C$ printed on the windscreen of your car,
the levels on concentric lines around some point at the bottom --where the
root of $T$ is situated and the windscreen wiper is fixed--
and if $s,t\in T\uhr C$ are of height $\alpha$ then $s$ stands to the left of
$t$ if and only if $s<_\alpha t$ according to $<_\alpha\in W$.
In this picture each position of a windscreen wiper,
which has its axis fixed in the root of $T$, corresponds to
a member $a$ of the generating chain $K_W$ induced by the wiper $W$.
And this member $a\in K_W$ can be calculated as the sum over the ``wiped
area'', i.e. the sum of all nodes which lie to the left of the windscreen
wiper. 

% The following lemma formally establishes this relationship between generating chains and wipers.
Establishing the relationship between wipers and generating chains
of a Souslin algebra is quite straight forward. Part b) of the next
lemma, which is a refinement of Lemma \ref{lm:genchaintree} requires a
second catch-up argument in order to achieve the denseness property (ii) of the
definition of wipers.
\begin{lm}[cf. {\cite[2.2.4]{diss}}]
\label{lm:genchainwiper}
Let $\B$ be a Souslin Algebra with Souslinisation $T$.
\begin{enumerate}[a)]
\item If $W$ is a wiper on $T$, then
$K_W$
%:=\big\{\sum M\mid(\exists\alpha<\omega_1)M\subseteq T_\alpha,\,
%M\text{ is an initial segment of }<_\alpha\big\}$$
is a generating chain of $\B$.
\item Let $K$ be a generating chain of $\B$.
Then there is a club $C\subseteq\omega_1$,
s.t. there is a wiper $W$ on $T\uhr C$,
inducing $K$, i.e. with $K=K_W$ as above.
\end{enumerate}
\end{lm}
By this correspondence we also see that
a small Souslin algebra always has exactly $2^{\aleph_1}$ maximal chains.

\section{Complete subalgebras and tree equivalence relations}

In order to represent complete, atomless subalgebras of s Souslin algebra we
will use certain equivalence relations defined on a Souslinisation.
The basic ideas for this method of representation can already be found in work
of Jech and Jensen from around 1970.
Here we will use the notions developed in \cite{diss}\footnote{Note that
  t.e.r.s as defined here and in \cite{SAE} are called \emph{decent}
  t.e.r.s in \cite{diss}. The properpty of decency of
  \cite[Def.1.1.3]{diss} refers to what is 
    called honesty in the present text and in \cite{SAE}.}
and \cite{SAE}
which are more appropriate for our purpose, that of rendering all Souslin
subalgebras isomorphic to the superalgebra.
In particular, we will take a topological point of view towards countable trees
which simplifies matters significantly.
(Cf. \cite[\S\S 2,5]{SAE} for a more detailed account on this.)

\begin{defi}
Let $T$ be a normal, $\aleph_0$-branching tree.
\begin{enumerate}[a)]
\item We say that an equivalence relation $\equiv$ on $T$
is a \emph{tree equivalence relation}
if $\equiv$ respects levels
(i.e., $\equiv$ refines $T\otimes T$),
is compatible with the tree order of $T$
(i.e., $\hgt s=\hgt r$ and $s<t\equiv u>r$ imply $s\equiv r$)
and honest (i.e. for each triple $s,s',t\in T$ such that $s<s'$ and $s\equiv
t$ we either find a node $t'>t$ such that $s'\equiv t'$ or there is no
immediate successor of $t$ equivalent to $s'\uhr(\hgt(s)+1)$, the immediate
successor of $s$ below $s'$) and
the quotient $T/\!\!\equiv$ is a normal
and $\aleph_0$-branching tree with the induced order.
\item A t.e.r. $\equiv$ on $T$ is said to be
\begin{itemize}
\item \emph{nice} if for all triples $s,s',t\in T$
with $s\equiv t$ and $s<s'$ there is some $t'>t$ equivalent to $s'$;
\item \emph{$\infty$-nice} if it is nice and for node $s$ and every
  $\equiv$-class $t/\!\!\equiv$ the intersection $T(s)\cap(t/\!\!\equiv)$ is
  either empty or infinite.(Equivavently we could require only
  $\NF(s)\cap(t/\!\!\equiv)$ to be infinite.)
\end{itemize}
\end{enumerate}
\end{defi}
A t.e.r. $\equiv$ on a Souslinisation $T$ of $\B$
represents the atomless and complete subalgebra
$\A_\equiv=\left\langle \sum t/\!\!\equiv\mid t\in T\right\rangle^\cm$ of $\B=\RO T$.
(We denote the $\equiv$-class of $t\in T$ by $t/\!\!\equiv$.) 
A subalgebra of $\B$ is called ($\infty$-)nice,
if it is represented by a(n $\infty$-)nice t.e.r. on some Souslinisation of
$\B$. 
Since niceness descends to the restriction to a club set of levels,
this is in\-dependent from the choice of the Souslinisation. Nice
  subalgebras are much easier to handle than general subalgebras. 
But as chain homogeneous Souslin algebras are necessarily homogeneous
(cf. Proposition \ref{prp:properties}) and 
every homogeneous, small Souslin algebra also has non-nice, atomless
and complete subalgebras
(\cite[Thms.3 \& 4]{SAE})
we also have to consider this general case in our constructions.
% Keep in mind, that \emph{we aim at rendering all atomless, complete
% subalgebras isomorphic to the Souslin algebra we intend to construct}.

Lemma 3 of \cite{SAE} states
that for every Souslin algebra $\B$ with a Souslinisation
$T$ and every complete and atomless subalgebra $\A$ or
$\B$ there is a t.e.r. $\equiv$ that represents $\A$
on $T\uhr C$ for some club $C$ of $\omega_1$.
If $\A$ is nice or $\infty$-nice then $\equiv$ can be chosen in a correponding
fashion.

A maximal branch $b$ of a tree $T$ is called \emph{cofinal} if
its order type with respect to $<_T$ coincides with the height of $T$.
We let
$$[T]=\left\{b\mid b\text{ is a cofinal branch of }T\right\}.$$
For $s\in T$ set $\hat{s}:=\{b\in[T]\mid s\in b\}$. We consider $[T]$
as a topological space with the base $\{\hat{s}\mid s\in T\}$.
Then for a normal tree $T$ of countable limit height $\alpha$ with only
countable levels the space $[T]$ is a perfect Polish space,
i.e., a separable and completely metrizable space without
isolated points.

It is clear that every tree homomorphism $\varphi:T\to S$
induces a continuous mapping $\overline{\varphi}:[T]\to[S]$ and that
for every club $c\subseteq\hgt T$ we have a natural homeomorphism
between $[T]$ and $[T\uhr c]$ since every cofinal branch of $T$
is uniquely determined by its intersection with $T\uhr c$.
If $T$ carries a t.e.r. $\equiv$, then an equivalence relation on $[T]$ is
induced which we also will denote by $\equiv$.

We now come to state the key property and the associated lemma we use
for constructions of
Souslin algebras with tightly controlled subalgebras.

% Note that, by an argument as in the proof of
% part a) of Proposition \ref{prp:decent_dense},
% given an $\infty$-nice t.e.r. $\equiv$ on $T\uhr\alpha$
% and any $x\in[T\uhr\alpha]$,
% then the $\equiv$-class of $x$ is
% a nowhere dense, perfect subset of the Polish space $[T\uhr\alpha]$
% and thus a perfect Polish space itself.

\begin{defi}
Let $X$ be a topological space and $\equiv$ an equivalence relation on $X$.
We say that a subset $N\subseteq X$ is \emph{suitable for} $\equiv$ if for
every element $x\in X$ the intersection of its $\equiv$-class $x/\!\!\equiv$
with $N$ is either empty or dense in $x/\!\!\equiv$ (with the subspace topology).
% For a mapping $h:X\to Z$ we say $N\subset X$ is \emph{suitable for }$h$ if $N$ is suitable
% for the equivalence relation given by $x\equiv y\, :\iff h(x)=h(y)$.
\end{defi}

In order to add a new limit level to a countable tree $T$ we will choose a
countable, dense subset $Q$ of $[T]$ and assign to
every member $x\in Q$ a bounding limit node.
A given a t.e.r. $\equiv$ will be extended to the new level
if and only if the set $Q$ is suitable for $\equiv$. Other wise it
would fail to be honest.
We will always succeed in extending $\equiv$ if we are allowed to
choose the members of $Q$ from a 
comeagre subset $M$ of $[T]$ that is suitable for $\equiv$.
The Reduction Lemma states that we will find such a set $M$.

\begin{lm}[{Reduction Lemma, \cite[Lemma 7]{SAE}}]\label{lm:SRL}
Let $T$ be a countable, normal and $\aleph_0$-branching
tree of limit height carrying a t.e.r. $\equiv$, and
% The induced equivalence relation on $[T]$ will also be called $\equiv$.
let $M\subseteq[T]$ be comeagre.
Then there is a comeagre subset $M'\subseteq M$ which is suitable for $\equiv$.
\end{lm}

For a monotone map $\varphi:S\to T$ between trees we denote the
mapping between the branch spaces $[S]\to[T],\,x\mapsto \varphi''x$ by
$\overline{\varphi}$. We will say that $\varphi$ induces $\overline{\varphi}$.
Extending a tree mapping $\varphi$ to a new limit level requires, in
the notation used above the Reduction Lemma, that the set $Q$ is
closed under the application of $\overline{\varphi}$.
The two following results establish that the Reduction Lemma can also
be applied to achieve this situation.

\begin{lm}
\label{lm:epimocomea}
For a normal, countable tree $T$ of limit height equipped with a t.e.r.
$\equiv$ consider the canonical tree epimorphism
$$\pi:T\to T/\!\!\equiv,\quad\pi(t)=t/\!\!\equiv.$$
% \begin{enumerate}[a)]\item 
Let the set
$Y:=\{y\in[T/\!\!\equiv]\mid (\exists x\in[T])\,y=\pi''x\}$
be equipped with the subspace topology inherited from $[T/\!\!\equiv]$.
Then the induced map
$$\overline{\pi}:[T]\to Y,\,x\mapsto\pi''x$$
is an open and surjective mapping and $Y$ is a comeagre subset of $[T]$.
% \item Let $\varphi:T/\!\!\equiv\,\to S$ be a tree isomorphism,
% $\rho:=\varphi\circ\pi$ and $Z:=\{\rho''x\mid x\in[T]\}$.
% Then
% $$\overline{\rho}:[T]\to Z, b\mapsto\rho''b$$
% is a continuous and open surjection, and $Z$ is comeagre in $[S]$.
% \end{enumerate}
\end{lm}
\begin{proof}
% As in the proof of Proposition \ref{prp:decent_dense},
Surjectivity of $\overline{\pi}$ as in the statement of the lemma is trivial.
In order to prove that $\overline{\pi}$ is an open mapping, we show that
the images of the basic open sets $\hat{t},\ t\in T$ are open in $Y$.
So let $t$ be a node of $T$. We prove the equation
$\overline{\pi}''\hat{t}=\bigcup\{\widehat{\pi(r)}\mid\,r\in \NF(t)\}\cap Y$.
Letting $\overline{\pi}(x)=\pi''x$ with $t\in x\in[T]$
be a member of the left hand side we choose $r$ to be the unique immediate
successor of $t$ in $x$ and see that this branch is also a member of the right
hand side.

For the converse inclusion let $y\in Y$ be in $\widehat{\pi(r)}$ for some
immediate successor $r$ of $t$. We have to find a branch $x\in[T]$ through $t$
with $\overline{\pi}(x)=y$. As $y$ is in $Y$ it has a $\overline{\pi}$-preimage
$x^*\in[T]$. Letting $\gamma:=\hgt(t)$, take $s$ to be the node of $x^*$ on
level $\gamma$ and $s'$ its immediate successor in $x^*$.
Then clearly $s\equiv t$ and $s'\equiv r$. Because of $\equiv$'s honesty,
there are on every level $\delta>\gamma$ successors of $t$ which stand
in relation $\equiv$ to a member of $x^*$. The set of all these nodes forms
(together with $t$ and its predecessors) a subtree $T^*$ of $T$ which
has a cofinal branch $x$ by K\"o{}nig's Lemma, because in any case we
can reduce $T^*$ to a club set of levels of order type  $\omega$
in which we find a finitely branching subtree of the same height.
Clearly we have $\overline{\pi}(x)=\overline{\pi}(x^*)=y$, as desired.

Finally, $Y$ is separable, metrizable and the image of the Polish space $[T]$
under the open mapping $\overline{\pi}$, so
$Y$ is Polish by a theorem of Sierpi\'{n}ski (cf. \cite[8.19]{kechris}) and
thus comeagre in $[S]$ by Choquet's Theorem (cf. \cite[8.17.ii]{kechris})
being a dense subset of $[S]$ that contains a Polish space as a subset.
\end{proof}

\begin{prp}\label{prp:comeagre_images}
Let $X$ and $Y$ be Polish spaces and $h:X\to Y$ a continuous mapping, such that
$h''X$ is comeagre in $Y$ and the right-hand side restriction of $h$ to its image,
i.e. $h:X\to h''X$ is an open mapping.
\begin{enumerate}[a)]
\item If $D\subset Y$ is meagre, then ${h^{-1}}''D$ is meagre as well.
\item If $M\subset X$ is comeagre, then $h''M$ is comeagre as well.
\end{enumerate}
\end{prp}
\begin{proof}
The proof of part a) is straight forward.
For part b) we can assume without loss of generality,
that $M$ is a dense $G_\delta$-subset of $X$.
Then clearly $h''M$ is an analytic subset of $Y$.
So, by a theorem of Lusin and Sierpinski
(cf.~\cite[Thm.11.18.b]{jechneu} or \cite[Thm.21.6]{kechris}),
$h''M$ has the Baire property,
i.e. there is an open subset $U$ of $Y$,
such that the symmetric difference of $h''X$ and $U$ is meagre.
In particular the sets $h''X\setminus U$ and $U\setminus h''X$ are meagre.
If we can show, that $U$ is dense in $Y$, then the proof is finished.

So assume to the contrary that there is
an open subset $V$ of $Y$ which is disjoint from $U$.
Then $V\cup h''X$ is meagre and so is
$M\cap {h^{-1}}''V={h^{-1}}''(V\cap h''M)$ by part a).
But since $M$ is comeagre and ${h^{-1}}''V$ is open, their
intersection cannot be meagre --- contradiction!
\end{proof}

We combine the Reduction Lemma and 
the last two propositions in order to get
the result which is appropriate for our use in the constructions.

% The proof is once again straight forward.

\begin{cor}
\label{cor:reduction}
Let $T$ be a countable, normal and $\aleph_0$-branching tree
of limit height $\alpha$.
Let $H$ be a countable set of triples $h=(c_h,\equiv_h,\varphi_h)$
where $c_h$ is a club subset of $\alpha$, $\equiv_h$ is a 
t.e.r. on $T\uhr c_h$ (with associated canonical mapping $\pi_h$) and
$\varphi_h:(T\uhr c_h)/\!\!\equiv_h\to T\uhr c_h$ is an isomorphism.
Let furthermore $I$ be a countable set of pairs $i=(c_i,\equiv_i)$
such that $c_i$ is club in $\alpha$ and $\equiv_i$ is an $\infty$-nice
t.e.r. on $T\uhr c_i$.
If $M$ is a comeagre subset of $[T]$ then there is a comeagre
subset $N$ of $M$, such that
\begin{enumerate}
\item  for all $h\in H$, $N$ is suitable for $\equiv_h$
and $(\overline{\varphi}_h\circ\overline{\pi}_h)''N=N$;
\item $N$ is suitable for all $\equiv_i$ where $i\in I$.
\end{enumerate}
\end{cor}
\begin{proof}[Proof idea]
As images and pre-images of comeagre sets under the mappings from $H$ are all
comeagre, repeated intersections deliver $N$.
\end{proof}

Our last technical lemma on t.e.r.s concerns nested $\infty$-nice t.e.r.s and
has a straight forward proof.
If $\equiv_1$ refines $\equiv_0$ as an equivalence relation on some set $M$,
we define the equivalence relation $\equiv_{0/1}$ on $M/\!\!\equiv_1$ by
$$(x/\!\!\equiv_1)\equiv_{0/1}(y/\!\!\equiv_1):\iff x\equiv_0 y$$
for $x,y\in M$.

\begin{lm}\label{lm:2ter}
Assume we are given two $\infty$-nice t.e.r.'s $\equiv_0$ and $\equiv_1$
on $T$ which is  normal, $\aleph_0$-branching and of
countable limit height,
such that $\equiv_0$ is refined by $\equiv_1$ and
the t.e.r. $\equiv_{0/1}$ on the quotient tree $T/\!\!\equiv_1$
induced by $\equiv_0$ is $\infty$-nice.
Let furthermore $M$ be a comeagre subset of $[T]$ that is suitable
for both $\equiv_0$ and $\equiv_1$.
Then $M/\!\!\equiv_1$ is suitable for $\equiv_{0/1}$.
\end{lm}

\section{Variations of Kurepa's Lemma}

A striking property of normal trees is formulated in
Kurepa's Isomorphism Lemma, cf.~\cite[p.102]{kurepa}.
The result and its proof are well-known,
but since we will use some variations of the argument later on,
we also sketch the proof.

\begin{lm}\label{lm:nor_tre_iso}
Let $S,\,T$ be two normal $\kappa$-branching trees, $\kappa\leq\omega$, of
the same limit height $\alpha<\omega_1$ with countable levels only.
Then $S\cong T$.
\end{lm}

\begin{proof}
% For trees of finite height, an easy inductive choice gives the isomorphism.
So let the height of $S$ and $T$ be a countable limit ordinal $\alpha$.
Choose countable and dense sets $X\subset[S]$ and $Y\subset[T]$
and enumerate them by
$(x_i\mid i\in\omega)$ and $(y_i\mid i\in\omega)$ respectively.

We give a back-and-forth-construction of a bijective mapping
$f:X\to Y$ which lifts to a tree isomorphism $\varphi:S\to T$.
Define $f(x_0)=y_0$ and
$\varphi(x_0\uhr\gamma)=y_0\uhr\gamma$
for all $\gamma<\alpha$.

As the forth-argument is completely analogous, we only explain hw the
back-step works.
Let $i$ be minimal such that $f^{-1}(y_i)$ has not yet been fixed, and
pick the minimal $\gamma$
such that $\varphi^{-1}(y_i\uhr\gamma)$ has not yet been defined.

This $\gamma$ is a successor ordinal, say $\gamma=\delta+1$,
because $S$ and $T$ are assumed to be normal trees.

Now choose an immediate successor $s$ of $\varphi^{-1}(y_i\uhr\delta)$ such that
$\varphi(s)$ has not yet been defined.
Such a node $s$ exists by the choice of $\gamma$.
Finally let $j$ be minimal such that $s\in x_j$.
Then $f(x_j)$ has not yet been defined,
but we set $f(x_j)=y_i$ and $\varphi(x_j\uhr\gamma)=y_i\uhr\gamma$.
This is consistent with the choices met for $f$ and $\varphi$ so far.

After every step of the construction $\varphi$ is a partial isomorphism
between $S$ and $T$. By the choice of $X$ and $Y$, the union of these
partial isomorphisms is bijective, so in the end $\varphi:S\to T$
is an isomorphism.
% If the trees are of infinite successor height $\alpha+n$, then find
% $\varphi:S\uhr\alpha\to T\uhr\alpha$ as above
% with the choices $X=\left\{\{r<_S s\}\mid s\in S_\alpha\right\}$
% and $Y=\left\{\{r<_T t\}\mid t\in T_\alpha\right\}$.
% Then $f$ gives the extension of $\varphi$ to level $\alpha$.
% The final $n$ steps are then easy.
\end{proof}

The argument also shows
that for any $\gamma<\alpha$ and every isomorphism
$\varphi_0$ between $S\uhr(\gamma+1)$ and $T\uhr(\gamma+1)$ there is an
extension $\varphi_0\subset\varphi:S\cong T$.

The denseness requirement in Defintion \ref{defi:wipers} of \emph{wipers}
enables us to combine the argument from the proof of Kurepa's Isomorphism
Lemma with that of Cantor's for the $\aleph_0$-categoricity of the dense
linear orders in order to construct isomorphisms between wipers.
% This will be a major step in the construction of a Souslin algebra that is homogeneous
% for all generating chains.

\begin{lm}\label{prp:extendwiperisos}
Let $T$ and $S$ be two countable, $\aleph_0$-branching, and normal
trees of the same height $\alpha<\omega_1$
and let $W_0=\langle\prec_\gamma\mid\gamma<\alpha\rangle$ be a wiper on $T$ and
$W_1=\langle <_\gamma\mid\gamma<\alpha\rangle$ a wiper on $S$.
Furthermore let $\beta<\alpha$ and $\varphi'$ be an isomorphism from
$T\uhr(\beta+1)$ onto $S\uhr(\beta+1)$, such that
for all $\gamma\leq\beta$ and $s,t\in T_\gamma$ we have
$$s\prec_\gamma t\Leftrightarrow \varphi'(s)<_\gamma\varphi'(t).$$
Then there is an isomorphism $\varphi$ between $T$ and $S$ extending $\varphi'$,
such that for all $\gamma<\alpha$ and $s,t\in T_\gamma$ we have
$$s\prec_\gamma t\Leftrightarrow \varphi'(s)<_\gamma\varphi'(t).$$
\end{lm}

\begin{proof}
We refer to the proof of Kurepa's Isomorphism Lemma
\ref{lm:nor_tre_iso}
and describe the only manipulation:
When it comes to choosing of $\varphi^{-1}(s)$,
this choice has to respect the wipers,
which is always possible by the denseness requirement (ii) of
Definition \ref{defi:wipers}.a).
\end{proof}

% In the construction of a big Souslin algebra in Section \ref{sec:big}
% we will carry out an iteration of length $\omega_2$.
% In the remainder of this section we prepare
% the treatment of certain subalgebras of the initial algebra of that iteration.

% Recall %from Section \ref{sec:nowhere_large} 
% that
% we call a t.e.r. $\equiv$ \emph{$\infty$-nice} if
% it is nice and for all $\alpha<\beta<\hgt(T)$ and all $s\in T_\beta$
% the projections $t\mapsto t\uhr\alpha$
% are $\infty$-to-one
% when restricted to the $\equiv$-class of $s$, i.e.,
% \begin{quote}
% for all $r\in (s\uhr\alpha)/\!\!\equiv$ the set 
% $\{t\in s/\!\!\equiv\,\mid t\uhr\alpha=r\}$ is infinite.
% \end{quote}
% A subalgebra $\A$ of $\B$ is called \emph{$\infty$-nice}
% if $\B$ has a Souslinisation $T$ which carries an $\infty$-nice t.e.r. $\equiv$
% such that $\{\sum s/\equiv\mid s\in T\}$ is dense in $\A$.

% The equivalence classes of $\infty$-nice t.e.r.s have an especially nice structure:
Our last Kurepa-style lemma states that isomorphisms between quotient trees by
$\infty$-nice t.e.r.s can be lifted to isomorphisms of the underlying trees.
In order to prove it, we first formulate a regularity property of subtrees
induced by an $\infty$-nice t.e.r.

\begin{prp}
Let $T$ be a normal, $\aleph_0$-branching tree of countable limit height
carrying an $\infty$-nice t.e.r. $\equiv$.
Then for every cofinal branch $b$ of $T$ the set
$\{s\in T\mid (\exists r\in b)s\equiv r\}$
is an $\aleph_0$-branching, normal tree.
\end{prp}
\begin{proof}
Unique limits are inherited from $T$, infinte branching follows from the
$\infty$-condition, successors in every higher level follow from niceness.
\end{proof}

% Recall Kurepa's Isomorphism Lemma \ref{lm:nor_tre_iso}
% which states that $\aleph_0$-branching normal trees of the
% same countable height are isomorphic.
% We once more refine Kurepa's Isomorphism Lemma in the following proposition,
% which eventually enables us to produce the additional
% isomorphisms we need in the construction in Section \ref{sec:big}.

\begin{prp}\label{lm:quo_tre_iso}
Let $T$ be a normal, $\aleph_0$-branching tree of limit height $\alpha<\omega_1$ and
let $\equiv$ and $\sim$ be $\infty$-nice t.e.r.'s on $T$. Let $\gamma$ be equal to $\alpha$ or else be a
successor ordinal below $\alpha$.
Let $\varphi'$ be an ismorphism between $(T\uhr\gamma)/\!\!\equiv$ and $(T\uhr\gamma)/\!\!\sim$.
Then there is an automorphism $\varphi$ of $T$ that carries $\equiv$ to $\sim$, s.t.
the induced map on $T/\!\!\equiv$ is an isomorphism onto $T/\!\!\sim$ that extends $\varphi'$.
\end{prp}
\begin{proof}
First use Lemma \ref{lm:nor_tre_iso} to extend $\varphi'$ to an isomorphism between
$T/\!\equiv$ and $T/\!\sim$ in the case that $\gamma<\alpha$. So we can assume $\alpha=\gamma$.

We now give a back-and-forth-argument
which lifts the isomorphism $\varphi'$
to an automorphism $\varphi$ of $T$.
Enumerate $T$ in order type $\omega$ by $s_0,s_1,\ldots$
For the induction step in the forth-direction
let $i$ be the minimal index $i$ for which
$\varphi(s_i)$ has not yet been defined.
Our choice for $\varphi(s_i)$
has to respect the tree order and the equivalence relations.
So let $s$ be the predecessor of maximal height
whose image under $\varphi$ is already determined by the choices met
earlier in the construction, where $s=s_i$ is allowed.
In the case that $s<_T s_i$,
we want to pick a node
$t\in \varphi'(s_i/\!\!\equiv)$ above $\varphi(s)$ which 
exists by the niceness of $\sim$.
We furthermore require that $t$
has not yet been assigned as some $\varphi(s_j)$.
This choice is possible due to the $\infty$-part in the
$\infty$-niceness of $\sim$.

For the back-step we replace in the above argument all $\varphi$ and
$\varphi'$ by $\varphi^{-1}$ and $\varphi'^{-1}$.
\end{proof}

\section{Statement of the main results}
\label{sec:hom_for_max}

\begin{defi}
Let $\B$ be a Souslin algebra. By \emph{continuous, $\infty$-nice
  chain of subalgebras of $\B$} we denote a sequence $(A_\nu\mid\nu<\beta)$ 
of $\infty$-nice subalgebras of $\B$ such that
\begin{enumerate}[(i)]
\item $\A_\mu$ is an $\infty$-nice-subalgebra of $\A_\nu$ for $\mu<\nu<\beta$,
\item for limit $\alpha<\beta$ the algebra $\A_\alpha$ is completely generated
  by the preceeding members $\A_\nu$, $\nu<\alpha$ of the chain, i.e.
  $$\A_\alpha=\left\langle\bigcup_{\nu<\alpha}\A_\nu\right\rangle^\cm.$$
\end{enumerate}
\end{defi}

\begin{thm}
\label{thm:hom_for_max}
Assume $\diaplus$. Then there is a small, chain homogeneous Souslin algebra
$\B$. Furthermore, the construction can be refined to achieve the following
additional homogeneity properties.
\begin{enumerate}[a)]
\item For every isomorphism $\varphi':\A\cong\A'$ of $\infty$-nice subalgebras
  of $\B$ there is an automorphism $\varphi$ of $\B$ extending $\varphi'$.
\item There is a continuous, $\infty$-nice chain $(\A_i\mid i<\omega_1)$ of 
  subalgebras of $\B$ such that:
  \begin{enumerate}[(i)]
  \item $\bigcup_{j<\omega_1}\A_j=\B$; 
  \item for every continuous, $\infty$-nice chain $(C_\nu\mid\nu<\lambda+1$)
    of $\B$ ($\lambda$ a countable limit) there is an isomorphism
    $\varphi:\A_\lambda\to\C_\lambda$ such for all $\nu<\lambda$ the
    restriction of $\varphi$ to $\A_\nu$ is an isomorphism onto $\C_\nu$.
%     for every normal sequence of countable ordinals
%     $(j_\nu\mid\nu<\lambda)$ with supremum $i$ there is an isomorphism
%     $\varphi:\A_\lambda\to\A_i)$ with $\varphi''\A_\nu=\A_{j_\nu}$. {\tt Hier
%       muss noch was passieren!}
  \end{enumerate}
\end{enumerate}
\end{thm}
The additional homogeneity properties properties stated in a) and b) of
Theorem \ref{thm:hom_for_max} have as 
primary goal their use in the construction which proves the following theorem.
It answers a question of Stevo Todor\c{c}evi\'{c}, who asked on the occasion
of a talk the author gave 
on chain homogeneous Souslin algebras at Toposym X in 2006 in Prague, whether
there are also big, chain homogeneous Souslin algebras.

\begin{thm}\label{thm:big}
Assume $\diaplus$ and let $\B$ be a chain homogeneous Souslin algebra as in
Theorem \ref{thm:hom_for_max}.
Then there is a big, chain homogeneous Souslin algebra $\B^*$ all of whose
subalgebras which are small Souslin algebras are isomorphic to $\B$.
\end{thm}
As the proof of this latter theorem is rather short, we give it right away.
\begin{proof}[Proof of Theorem \ref{thm:big} from Theorem
  \ref{thm:hom_for_max}]
We realise the big and chain homogeneous Souslin algebra $\B^*$
as the union of an increasing chain of small,
chain homogeneous Souslin algebras $(\B_\alpha\mid\alpha<\omega_2)$.
All the small algebras $\B_\alpha$ on this chain are isomorphic to the initial
algebra $\B_0:=\B$, which we take from Theorem \ref{thm:hom_for_max}.
The properties of $\B^*$ as stated in the theorem then follow from
general principles. As all the $\B_\alpha$ satisfy the c.c.c. and $\B^*$ is
defined as their direct limit, $\B^*$ also satisfies the c.c.c. and is therefore
complete. Distributivity and chain homogeneity are also inherited.

As above, denote by $(\A_i\mid i<\omega_1)$ the continuously increasing chain 
of $\infty$-nice subalgebras of $\B$.
For the successor step, if $\B_\alpha$ is given,
choose any isomorphism $\varphi_\alpha:\B_\alpha\cong\A_0$.
Next choose a Souslin algebra
$\B_{\alpha+1}$ extending $\B_\alpha$,
such that there is an extension $\varphi_{\alpha+1}$
of $\varphi_\alpha$ and witnessing $\B_{\alpha+1}\cong\B$.
This extension exists by part {\it (a)}  of Theorem \ref{thm:hom_for_max}.

If $\lambda<\omega_2$ is of countable cofinality,
we choose a normal sequence $(i_\nu\mid\nu<\mu)$
with $\sup_{\nu<\mu} i_\nu=\lambda$
for some countable limit ordinal $\mu$.
Inductively construct an increasing chain $(\psi_\nu\mid\nu<\mu)$
of isomorphisms $\psi_\nu:\A_\nu\cong\B_{i_\nu}$.
Then choose $\B_\lambda$ as a super-algebra of
$$\bigcup_{\nu<\mu}\B_{i_\nu}=\bigcup_{\alpha<\lambda}\B_\alpha$$
isomorphic to $\A_\mu$ by some extension $\psi_\mu$ of
$\bigcup_{\nu<\mu}\psi_\nu$. 

If $\lambda$ has uncountable cofinality we simply let
$\B_\lambda=\bigcup_{\alpha<\lambda}\B_\alpha$, which in this case is the
direct limit. 
Then we see that $\B_\lambda\cong\B_0$ for $\lambda<\omega_2$
by choosing a cofinal sequence $(i_\nu\mid\nu<\omega_1)$
and recursively choosing a chain of isomorphisms
$\psi_\nu:\A_\nu\cong\B_{i_\nu}$.
This goes through limit stages $\nu<\omega_1$ by property {\it(b.ii)}
of Theorem \ref{thm:hom_for_max}.
% The direct limit $\B_{\omega_2}$ of the increasing chain
% of small Souslin algebras is Souslin as well:
% the c.c.c. is preserved by direct limits
% (i.e. finite support iterations or here simply the union of the $\B_\alpha$).
% And then the same holds for distributivity,
% for if $(a_{ij}\mid i,j<\omega)$ is any family of members of $\B_{\omega_2}$
% then there is some $\alpha<\omega_2$
% with $a_{ij}\in\B_\alpha$ for all $i,j<\omega$,
% and $\B_\alpha$ witnesses that the distributivity law is preserved.
Iterating up to $\alpha=\omega_2$ finishes the construction of $\B^*$.
\end{proof}

\section{Final preparations}
\label{sec:finalprep}

% diaplus-machinery
In our Souslin tree construction
we want to build additional objects (mappings, tree isomorphisms)
on club sets of levels of the tree $T$ to be constructed,
that relate given objects (e.g. pairs of wipers or t.e.r.s on $T\uhr C$) to
each other. 
During a relevant construction step,
initial segments of some of the given objects
are proposed by a $\diaplus$-sequence
and in order to extend an additional object we need some pointer
to indicate the ordinal stage up to which the recursive construction
of the additional object has reached so far. This is the function
of $\varepsilon$ of the following definition.
\begin{defi}[the $\diaplus$-machinery]
Fix a $\diaplus$-sequence $(S_\alpha)_{\alpha<\omega_1}$.
\begin{enumerate}[a)]
\item Let $C_0:=\{\alpha<\omega_1\mid \omega\alpha=\alpha\}$ be the set of countable fixed points
of the left-multiplication with $\omega$.
\item Let $\alpha\in C_0$. For $x\in\pome(\alpha)$ set
$$c(x)=\{\gamma\in C_0\cap\alpha\mid
x\cap\omega(\gamma+1)\setminus\omega\gamma\neq\varnothing\}.$$
(This will the set of levels ofthe tree, to which $x$ refers.)
\item The \emph{set of relevant guesses for stage} $\alpha$, $G_\alpha$
is the set of pairs
$(x,d)\in S_\alpha\times S_\alpha$ such that
$c(x),\,d$ and $c(x)\cap d$ are club in $\alpha$
and for $\gamma\in d$ the sets $x\cap\gamma$
and $d\cap\gamma$ are in $S_\gamma$.
\item For $(x,d)\in G_\alpha$ let
$$e_{x,d}:=\{\gamma\in c(x)\cap d\mid\bigcup\left(\gamma\cap c(x)\cap d\right)=\gamma\}$$
be the Cantor-Bendixson derivative of $c(x)\cap d$, i.e., the set of its limit points.
\item Let $\varepsilon_{x,d}:=\bigcup e_{x,d}$. (Note that $\varepsilon_{x,d}=0$ if $\ot(c(x)\cap d)=\omega$.)
\end{enumerate}
\end{defi}
Now if for example $(x,d)\in G_\alpha$ and $x$ codes a pair of wipers on $T\uhr c(x)$,
we know that up to stage $\varepsilon_{x,d}$ our recursive construction of the additional
object --- here: an isomorphism between the wipers given by $x$ --- has been invoked and
therefore up to this stage this isomorphism has yet been constructed.

% Engaging relations
Recall Lemma \ref{lm:quo_tre_iso},
which states that an isomorphism between two quotient trees
$T/\!\!\equiv_0$ and $T/\!\!\equiv_1$ can be lifted to
an automorphism of $T$ if the t.e.r.s are $\infty$-nice.
In order to get hold of isomorphisms between quotient trees
we introduce the notion of an engaging relation. 
Consider two $\infty$-nice t.e.r.'s $\equiv_0$ and $\equiv_1$ on
$T\uhr c$, with $c\subset\alpha$ club,
and an isomorphism
$\varphi:(T\uhr c)/\!\!\equiv_0\,\to(T\uhr c)/\!\!\equiv_1$
between the quotient trees.
Then $\varphi$ naturally induces a relation $\Phi$ on $T\uhr c$,
that consists of the pairs $s,t\in T\uhr c$ with
$\varphi(s/\!\equiv_0)=t/\!\equiv_1$.
The properties of such a relation $\Phi$ are captured
in the following definition.

\begin{defi}\label{defi:engaging}
We say that a relation $\Phi$ on a tree $T$ satisfying points 1-4) 
below is \emph{engaging}.
\begin{enumerate}
\item There is a set $c_\Phi$ such that for all $s,t\in T$
we have $s\Phi t$ only if $\hgt(s)=\hgt(t)\in c_\Phi$,
\item the left-induced relation $\Phi^0:=\{(s,s')\in(T\uhr c)^2\mid(\exists t)s\Phi t\text{ and }s'\Phi t\}$
is an $\infty$-nice t.e.r. on $T\uhr c_\Phi$,
\item the right-induced relation
$\Phi^1:=\{(t,t')\in(T\uhr c)^2\mid(\exists s)s\Phi t\text{ and }s\Phi t'\}$ also is an $\infty$-nice t.e.r.
on $c_\Phi$
\item $\Phi$ induces an isomorphism $\varphi_\Phi$ between $(T\uhr c)/\Phi^0$ and $(T\uhr c)/\Phi^1$ via
$\varphi_\Phi(s/\Phi^0)=t/\Phi^1$ for any $t$ with $s\Phi t$.
\end{enumerate}
\end{defi}
It is clear that the relation $\Phi$ considered above the last definition is indeed engaging with $c_\Phi=c$,
and $\Phi^0$ is $\equiv_0$ while $\Phi^1$ is $\equiv_1$.

% C_i family
Finally we define a certain sequence of club sets of $\omega_1$
which will act as the familiy of supporting sets of levels for the t.e.r.s
that represent the increasing sequence $(\A_i\mid i<\omega_1)$ if
$\infty$-nice subalgebras.
Recall the definition of $C_0$, the set of infinite
fixed points of the left hand ordinal multiplication with $\omega$ in $\omega_1$:
$$C_0:=\{\alpha<\omega_1\mid \alpha\neq 0,\,\omega\alpha=\alpha\}.$$
Inductively define $C_{i+1}$ for $i<\omega_1$ to be the Cantor-Bendixson-derivative of $C_i$,
and for limit ordinals $i$ let $C_i$ be the intersection of the $C_j$ defined so far:
$$C_{i+1}:=\{\alpha\in C_i\mid \sup(C_i\cap\alpha)=\alpha\}\quad\text{ and }\quad
C_i=\bigcap_{j<i}C_j\text{ for limit }i.$$
We list some properties of the sequence $(C_i)_{i<\omega_1}$ used in the construction below.
\begin{enumerate}[(a)]
\item all the $C_i$ are club in $\omega_1$,
\item the sequence is continuously decreasing and has an empty intersection,
hence there is for every $\alpha\in C_0$ a unique $i=i(\alpha)$
with $\alpha\in C_i\setminus C_{i+1}$,
\item every $\alpha\in C_0$ has a direct predecessor in $C_{i(\alpha)}$, call it $\alpha^-$,
\item for limit $i<\omega_1$, the minimum of $C_i$ is the supremum of the minima of the
$C_j$ for $j<i$.
\end{enumerate}

\section{Proof of Theorem \ref{thm:hom_for_max}}
\label{sec:mainproof}
Before we step into the details of the construction, we give a brief sketch of
the main steps to be carried out at one stage of the recursion.
We use the principle $\diaplus$ during the construction in order to get hold
of all kinds of objects we are interested in, which are
\begin{itemize}
\item pairs of wipers that represent pairs of generating chains;
\item t.e.r.s that represent Souslin subalgebras;
\item engaging relations that represent isomorphism between
  $\infty$-nice subalgebras.
\end{itemize}
All of these are relations on $T\uhr C$ for some club $C$ of $\omega_1$.
In a relevant stage of the construction, say stage $\alpha\in C_0$,
we have the tree $T\uhr\alpha$ at hand and carry out five steps:
\begin{enumerate}
\item We collect all pairs $(x,d)\in G_\alpha$ that correspond to any of the
  above items. If necessary, we extend the mappings associated to these
  objects which were chosen in earlier stages of the construction by virtue of
  Kurepa's Lemma and its variations up to below $T_\alpha$.
\item We also collect the relvant t.e.r.s which represent the members of our
  increasing sequence of $\infty$-nice subalgebras.
\item Now we consider the induced relations and mappings on the Polish space
  $[\Talpha]$ and find with the aid of the Reduction Lemma a comeagre subset
  $N\subset[\Talpha]$ that is suitable for all relevant equivalence relations
  and is equal to its image under all relevant mappings:
  $\overline{\rho}''N=N$.
\item Finally we choose a countable, dense subset $Q$ of $N$ with the same
  properties correponding to the objects of interest: suitable for the
  equivalence relations and $\overline{\rho}''Q=Q$. The members of $Q$
  represent the nodes in our new limit level $T_\alpha$ of the tree $T$.
\item If necessary we extend the t.e.r. which represents to relevant member of
  the $\infty$-nice chain.
\end{enumerate}

\subsection{Setting the stage, casting and styling}
We now begin and inductively, i.e., level by level, construct a
tree-order $<_T$ on the set $\omega_1$, such that
the resulting tree $T$ will be a normal and $\aleph_0$-branching Souslin tree.
The Souslin algebra to be constructed will be the regular open algebra of
$T=(\omega_1,<_T)$.

Let $0$ be the root of $T$ and in every successor step
fix $\aleph_0$ distinct direct successors for
each maximal node in such a way that
for every $\alpha<\omega_1\setminus\{0\}$ the level
$T_{\alpha}$ consists of the next $\omega$ many ordinals not yet
used in the construction. So we have
$T_1=\omega\setminus\{0\}$, 
$T_{n}=\omega n\setminus\omega(n-1)$
for all natural numbers $n\geq 2$ and finally $\Talpha=\omega\alpha$
for all infinite, countable ordinals $\alpha$.

We fix a $\dia$-sequence $(R_\nu)_{\nu<\omega_1}$, a $\diaplus$-sequence
$(S_\nu)_{\nu<\omega_1}$ and a bijection 
$g:\omega_1\rightarrow(2{}\times{}\omega_1{}\times{}\omega_1)$
with $g''\lambda=2{}\times{}\lambda{}\times{}\lambda$ for all limit ordinals
$\lambda$.
Let for $i\in 2$ be $g_i$ the concatenation of $g$ and the projection onto
the fibre over $i$.
We will use $g$ for coding t.e.r.s and pairs of wipers as sets
of ordinals.

Now let $\alpha<\omega_1$ be a limit ordinal.
By the inductive assumption we have so far
constructed a normal tree order $(T{}\upharpoonright{}\alpha,<_T)$ on the
supporting set $\omega\alpha$.

% We consider the space $[T{}\uhr{}\alpha]$
% of all cofinal branches of $T{}\uhr{}\alpha$
% with the topology generated by the clopen sets
% $\hat{t}:=\{z\in [\Talpha]:t\in z\}$ for $t\in T{}\uhr{}\alpha$.
% By Proposition \ref{prp:cofbrPolish} we know
% that this topology is Polish and perfect.

If $\alpha<\omega\alpha$ we simply choose a countable dense subset $Q_\alpha$ of
$[\Talpha]$ and embed $Q_\alpha$ onto $\omega(\alpha+1)\setminus\omega\alpha$,
i.e., we choose a bijection between $Q_\alpha$ and $\omega(\alpha+1)\setminus\omega\alpha$
and extend $<_T$ on $\omega(\alpha+1)$ in the obvious way.

If $\alpha=\omega\alpha$, we want to choose a countable dense subset
$Q_\alpha$ of $[\Talpha]$, too, but this time our set also has to seal certain
maximal antichains of $\Talpha$ and extend enough tree isomorphisms.
% In order to state the inductive assumption and
% the inductive claim we need some notation.
Recall the definitions from the last section concerning the use
of the $\diaplus$-sequence ($G_\alpha,\, c(x),\, \varepsilon_{x,d}$, etc.).

From now on, if $(x,d)\in G_\alpha$ is fixed, we write $c$ for $c(x)$.
Denote by $E_\alpha$ the subset  of $S_\alpha\times S_\alpha$ that consists of
all pairs $(x,d)\in G_\alpha$ such that
$g_0''x$ and $g_1''x$ are wipers
$W_0=\langle \prec_\gamma:\gamma\in c(x)\rangle$
and $W_1=\langle <_\gamma:\gamma\in c(x)\rangle$ on $T\uhr c(x)$ respectively.
% (It would have been more correct to write ${g_0}''x=\bigcup W_0$ and
% ${g_1}''x=\bigcup W_1$.) 
So $(x,d)\in E_\alpha$ if $x$ codes a pair of wipers on $\Tc$ and is
guessed correctly by the $\diaplus$-sequence along with a club $d$ on the
members of $d$ itself.

Now we take care of the homogeneity property stated above as {\it (a)}.
We define 
$$E'_\alpha:=\left\{(x,d)\in G_\alpha\mid x\text{ codes an engaging relation
  }\Phi_x\text{ on }T\uhr c(x)\right\}.$$ 
Every $(x,d)\in E'_\alpha$ induces a pair of $\infty$-nice t.e.r.s
$\equiv^x_0$ and $\equiv^x_1$ and an isomorphism $\varphi'_x$ between $(T\uhr
c)/\!\!\equiv^x_0$ and $(T\uhr c)/\!\!\equiv_1$.

The guesses of the $\diaplus$-sequence for t.e.r.s
are collected in the set $F_\alpha$:
$$F_\alpha:=\{(x,d)\in G_\alpha\mid x\text{ codes a t.e.r. on }c(x)\}.$$
Here we code e.g. with respect to the map $g_0$.
For $(x,d)\in F_\alpha$ let $\equiv_x$ denote the t.e.r. which is coded
by $x$.
It induces the $\aleph_0$-branching, normal tree $(T\uhr c)/\!\!\equiv_x$.

Keep in mind that $c\cap d$ is unbounded in $\alpha$ if $(x,d)\in E_\alpha$ or
$F_\alpha$ by the definition of the relevant stages set $G_\alpha$.

The inductive hypothesis is that
$\Talpha$ is a countable and $\aleph_0$-branching tree of height $\alpha$ and
\begin{enumerate}
\item for every pair $(x,d)\in E_\alpha$ there is a $\subset$-chain
  $\langle\varphi_{x\cap\gamma,d\cap\gamma}:\gamma\in e_{x,d}\rangle$ and each
  $\varphi_{x\cap\gamma,d\cap\gamma}$ is a tree automorphism of $T\uhr(c\cap
  d\cap\gamma+1)$ that was fixed in induction step $\gamma\in e_{x,d}$ and
  which carries $W^0_{x,d}\uhr c\cap d\cap\gamma+1=\langle
  \prec_\delta:\delta\in c\cap d\cap\gamma+1\rangle$ to $W^1_{x,d}\uhr c\cap
  d\cap\gamma+1=\langle <_\delta:\delta\in c\cap d\cap\gamma+1\rangle$. In
  short: $\varphi_{x\cap\gamma,d\cap\gamma}$ is an isomorphism of wipers between
  $W^0_{x,d}\uhr c\cap d\cap\gamma+1$ and $W^1_{x,d}\uhr c\cap d\cap\gamma+1$
  for all $\gamma\in e_{x,d}$;
\item for every pair $(x,d)\in E'_\alpha$ there is an automorphism of $T\uhr
  c\cap d\cap \varepsilon_{x,d}$ which extends $\varphi'_x$;
\item for every pair $(x,d)\in F_\alpha$ there is a $\subset$-chain
  $\langle\psi_{x\cap\gamma,d\cap\gamma}:\gamma\in e_{x,d}\rangle$ of tree
  isomorphisms $\psi_{x\cap\gamma,d\cap\gamma}:T\uhr(c\cap
  d\cap\gamma+1)/\!\!\equiv_x \to T\uhr(c\cap d\cap\gamma+1)$ each of them
  fixed in induction step $\gamma\in e_{x,d}$.
\end{enumerate}
(Recall that $c=c(x)$ and $e_{x,d}$ is the set of limit points of $c\cap d$.)
For those $(x,d)$ in $E_\alpha, E'_\alpha$ or in $F_\alpha$ but with
$\varepsilon_{x,d}<\alpha$ we need to choose extensions for the maps granted
by the inductive hypothesis as follows.

Fix $(x,d)\in E_\alpha$. If $\varepsilon_{x,d}=\alpha$ let
$$\varphi_{x,d}:=\bigcup_{\gamma\in
  e_{x,d}}\varphi_{x\cap\gamma,d\cap\gamma}.$$ 
Otherwise extend $\varphi_{x\cap\varepsilon_{x,d},d\cap\varepsilon_{x,d}}$
by Proposition \ref{prp:extendwiperisos} to some isomorphism $\varphi_{x,d}$
between $W^0_{x,d}\uhr c\cap d$ and $W^1_{x,d}\uhr c\cap d$.

% E' extension
% and (in the notation of Definition \ref{defi:engaging})
For $(x,d)\in E'_\alpha$ use Lemma \ref{lm:quo_tre_iso} to lift
$\varphi_{\Phi_x}$ to an automorphism $\varphi_{x,d}$ of $T\uhr c(x)\cap d$
carrying $\Phi^0$ to $\Phi^1$ which extends $\varphi'_x$.

Let $(x,d)\in F_\alpha$ and set
$\psi_{x,d}:=\bigcup_\gamma \psi_{x\cap\gamma,d\cap\gamma}$
if $\varepsilon_{x,d}=\alpha$.
Otherwise extend the union of the chain
by Kurepa's Lemma \ref{lm:nor_tre_iso} to some isomorphism
$$\psi_{x,d}:(T\uhr c\cap d)/\!\equiv_x\to T\uhr c\cap d.$$

\subsection{Arranging the scene for the $\infty$-nice chain}

We now describe how to embed the increasing $\omega_1$-sequence of
$\infty$-nice subalgebras $\A_i$ in $\B$.
Recall the definition of the decreasing sequence
$(C_i\mid i\in\omega_1)$ of club subsets of $\omega_1$
as well as the derived definitions of $i(\alpha)$ and $\alpha^-$.
The $\infty$-nice t.e.r. $\equiv_i$ representing $\A_i$
will be defined on $T\uhr C_i\cup\{0\}$ in the course
of the construction of $T$. 
The main requirements to meet are:
\begin{enumerate}
\item $\equiv_i$ is $\infty$-nice for all $i<\omega_1$,
\item for $i>j$ the restriction of $\equiv_j$ to $T\uhr C_i$
is refined by $\equiv_i$ in a way such that
the induced t.e.r. $\equiv_{j/i}$ on the normal Souslin tree
$(T\uhr C_i)/\!\equiv_i$ is $\infty$-nice,
\item for limit $i<\omega_1$ we want to have $\A_i=\langle\bigcup_{j<i}\A_j\rangle^\cm$, so
$\equiv_i$ shall be the conjunction of the $\equiv_j$ for $j<i$ in this case:
$$s\equiv_i t:\Leftrightarrow (\forall j<i)s\equiv_j t.$$
\end{enumerate}
For any $i$, on level $T_0=\{\text{\bf root}\}$
the relation $\equiv_i$ is of course trivial.
On level $T_{\min C_i}$ we define $\equiv_i$
to be the identity, i.e.,
$s\equiv_i t$ if and only if $s=t$ for $s,t\in T_{\min C_i}$.
This is a minor violation of the $\infty$-niceness requirement
we posed on $\equiv_i$.
But this is easily remedied by deleting $\min C_i$
from the club set $C_i$.
On the other hand,
by this convention we directly see that in the end
$\bigcup \A_i$ will be a dense subset of $\B$,
because $\{\min C_i\mid i<\omega_1\}$ is unbounded in $\omega_1$.

In level $\alpha\in C_0$ we have that for all $j<i:=i(\alpha)$
the set $C_j\cap\alpha$ is club in $\alpha$ and
the t.e.r. $\equiv_j$ on $T\uhr(C_j\cap\alpha)$
has by normality of the quotient tree
a unique t.e.r.-extension to $T_\alpha$.
So, to satisfy the niceness condition for the t.e.r.s $\equiv_j$ with $j<i(\alpha)$, 
level $T_\alpha$ has to be chosen suitably %carefully
with the aid of the Reduction Lemma \ref{lm:SRL}.
% In the case where $i(\alpha)$ is a limit ordinal,
% we even have $\equiv_{i(\alpha)}$ on $[\Talpha]$ at hand
% before we choose $T_\alpha$,
% by the requirement (3) above.
% It is easily seen that, since for a limit ordinal $i$
% the minimum of $C_i$ is just the supremum of the minima of the $C_j$ for $j<i$,
% this is consistent with our appointment
% that $\equiv_i$ be the identity on $T_{\min C_i}$.

Let now $\alpha\in C_0$ be such that $i=i(\alpha)$ is a limit ordinal. For $j<i$ define $\equiv_j$ on $[\Talpha]$
as above and let for $x,y\in [\Talpha]$
$$x\equiv_{i} y:\iff (\forall j<i) x\equiv_j y.$$
This is the coarsest possibility to extend $\equiv_{i}$ to $T_\alpha$ and the only one by our requirement
that $\bigcup_{j<i}\A_j$ completely generate $\A_{i}$.

We now show that the relation $\equiv_{i}$ on $[\Talpha]$ is induced by a
(single) $\infty$-nice t.e.r. $\simeq$. Therefore the Reduction Lemma
\ref{lm:SRL} can be applied to $\equiv_i$. The t.e.r. $\simeq$ will be defined
as a diagonal along the $\equiv_j$ for $j<i$.

Fix $\delta_0$ in $C_0$ between $\alpha^-$ and $\alpha$, so
$j_0=i(\delta_0)<i$. Then define 
$$j_\nu:=i(\delta_\nu)\quad\text{ and }\quad
\delta_{\nu+1}:=\min C_{j_\nu+1}\setminus\delta_\nu=
\min C_{j_\nu+1}\setminus\delta_0$$
and $\delta_\mu:=\sup\{\delta_\nu\mid\nu<\mu\}$ for limit ordinals $\mu$ with
$i(\delta_\nu)<i$ for all $\nu<\mu$. Then the $i(\nu)$ are the ordinals from
$j_0$ up to $i$. The final $\delta_\mu$ is just $\alpha$, and this ordinal
$\mu=\ot(i\setminus j_0)$ is a limit.

The set of the $\delta_\nu$ joined with $\{0\}$
is the club set of $\alpha$ on which we now define the t.e.r. $\simeq$.
For $s,t\in T_{\delta_\nu}$ define
$$s\simeq t:\iff [(\forall \beta<\delta_0)(i(\beta)\leq i\Rightarrow
s\uhr\beta\equiv_{i(\beta)}t\uhr\beta)\quad\text{and}\quad s\equiv_{j_\nu}t\,\,]$$ 
Then $\simeq$ is an $\infty$-nice t.e.r.: fix $s\simeq t$ on level $\delta_\xi$ and consider $s'>s$, where
$s\in T_{\delta_\nu}$.
In the successor case, letting $\nu=\nu^-+1$, i.e. $j_\nu=j_{\nu^-}+1$,
this follows from the definition of $\equiv_{j_\nu}$
having only classes that are dense subsets of $\equiv_{j_{\nu^-}}$-classes.
To argue for niceness in the limit case we refer to the choice of $T_{\delta_\nu}$ which assures,
that the $\equiv_{j_\nu}$-classes lie densely in the $\equiv_{j_\nu}$-classes while the $\infty$-part
of $\infty$-niceness is trivially satisfied on limit stages when satisfied everywhere below.
Since the $j_\nu$ are cofinal in $i$,
the $\infty$-nice t.e.r. $\simeq$ induces $\equiv_i$
on $[\Talpha]$.

However, for successor $i(\alpha)$
the definition of $\equiv_{i(\alpha)}$ on $T_{\alpha}$
involves the choice of the new level $T_\alpha$.
So we continue by giving the rules
for the choice of $T_\alpha$. % in the successor case $i=i(\alpha)=j+1$.
 
\subsection{Suitably reducing the area
%Reduction to a suitable comeagre and finally countable, dense set
}
Now Lemma \ref{lm:epimocomea} comes into play. Let for $(x,d)\in F_\alpha$
$$\pi_{x,d}:T\uhr c\cap d\to (T\uhr c\cap d)/\equiv_x,\quad s\mapsto s/\equiv_x$$
be the canonical mapping associated to $\equiv_x$ and define
$$\rho_{x,d}:=\psi_{x,d}\circ\pi_{x,d}:T\uhr c\cap d\to T\uhr c\cap d.$$
Since $\psi_{x,d}$ is an isomorphism,
% $\rho_{x,d}$ has the same properties as the map $\varphi$ in
% Lemma \ref{lm:epimocomea}. 
the induced continuous map
$$\overline{\rho}_{x,d}:[\Talpha]\to[\Talpha],\quad
b\mapsto\{s\mid (\exists t\in b\cap(T\uhr c\cap d))s<_T \rho_{x,d}(t)\}$$
has a comeagre image in $[\Talpha]$ and
is an open mapping when the range
is restricted to ${\overline{\rho}_{x,d}}''[\Talpha]$.

The $\dia$-sequence $(R_\nu)_{\nu<\omega_1}$ proposes candidates for maximal
antichains in the usual way.
If $R_\alpha$ is a maximal antichain of $\Talpha$ then we have to
ensure that each member of $T_\alpha$ is a $<_T$-successor of some
element of $R_\alpha$. That means $Q_\alpha$ has to be a subset of
$$M_\alpha=\{x\in [\Talpha]: \exists\gamma<\alpha\, x\uhr\gamma\in R_\alpha\}$$
which is itself an open dense subset of $[\Talpha]$, because $R_\alpha$ is a maximal
antichain.
If $R_\alpha$ is not a maximal antichain in $\Talpha$ we simply set $M_\alpha=[\Talpha]$.
We apply Corollary \ref{cor:reduction} of the Reduction Lemma to the sets
$M=M_\alpha$ and
$$H=\{\overline{\varphi}_{x,d}\mid (x,d)\in E_\alpha\}
\cup\{\overline{\rho}_{x,d}\mid (x,d)\in F_\alpha\}$$
and
$$I:=\begin{cases}
\{\equiv_j\mid j<i(\alpha)\}, & \text{ if }i(\alpha)\text{ is a successor
  ordinal}\\
\{\equiv_j\mid j\le i(\alpha)\}, & \text{ if }i(\alpha)\text{ is a limit
  ordinal}\\
\end{cases}.$$
Since the $E_\alpha\cup E'_\alpha$-part of $H$ consists of homeomorphisms and the
$F_\alpha$-part is subject to Lemma \ref{lm:epimocomea}, the hypotheses are
satisfied. The result is a comeagre subset 
$N_\alpha$ of $M_\alpha$, such that 
\begin{enumerate}[(1)]
\item for all $(x,d)\in E_\alpha\cup E'_\alpha$ we have
  $\overline{\varphi}_{x,d}''N_\alpha=N_\alpha$;
\item for each $(x,d)\in F_\alpha$, $N_\alpha$ is suitable for
  $\equiv_x$ and $\overline{\rho}_{x,d}''N_\alpha=N_\alpha$;
\item $N_\alpha$ is suitable for all $\equiv_i\in I$.
\end{enumerate}
% (Recall that we say that a comeagre subset $M$ of $[T\uhr\alpha]$
% is suitable for $\overline{\rho}_{x,d}$ if for every
% branch $b\in[T]$ the intersection of the class $b/\!\!\equiv_x$ with $M$ is
% either empty or dense in $b/\!\!\equiv_x$.)

% Check: ist das wirklich ueberfluessig?:
% Let $F_\alpha^*$ be the set of those $(x,d)\in F_\alpha$,
% such that $N_\alpha$ is suitable for $\equiv_x$.
% So we can choose for every $(x,d)\in F^*_\alpha$
% a right inverse $\sigma_{x,d}:N_\alpha\to N_\alpha$ of
% $\overline{\rho}_{x,d}\uhr N_\alpha$.
% In general, these sections $\sigma_{x,d}$ will not be continuous,
% but this is no longer important.
% To rule out the bad t.e.r.s,
% pick for each $(x,d)\in F_\alpha\setminus F^*_\alpha$
% one witness $b_x\in N_\alpha$ for the fact that
% $N_\alpha$ is not suitable for $\equiv_x$, i.e., that
% $N_\alpha\cap b_x/\!\!\equiv_x$ is not dense in $b_x/\!\!\equiv_x$.
% As already noted in Remark \ref{rem:FRL}, this choice
% impeaches every t.e.r. extending $\equiv_x$
% from being a \emph{decent} t.e.r.,
% because it fails to satisfy the necessary denseness condition stated
% in Proposition \ref{prp:decent_dense}.

The inductive claim 
is that there is a choice for $Q_\alpha$,
i.e., a countable and  dense subset of $M_\alpha$ that shares properties
(1-3) of $N_\alpha$ stated above. For then, $\Talpha$ can be extended by a
level $T_\alpha$ of nodes corresponding to the branches in $Q_\alpha$,
in a way that guarantees, that the relevant t.e.r.s and tree homomorphisms
extend to the new level.

% \begin{enumerate}
% \item every tree automorphism $\varphi_{x,d}$ of $\Tcd$ extends to $Q_\alpha$
%   for $(x,d)\in E_\alpha\cap E'_\alpha$ and
% \item for every $(x,d)\in F_\alpha$, the tree isomorphism
%   $\psi_{x,d}:(T\uhr c\cap d)/\!\!\equiv_x\to T\uhr c\cap d$ extends to the
%   respective trees with the new top level $T_\alpha$ corresponding to
%   $Q_\alpha$ (or $T_\alpha/\!\!\equiv_x$ respectively) added on.
% \item $Q_\alpha$ is suitable for all t.e.r.s $\$
% \end{enumerate}

First choose a left-inverse mapping $\sigma_{x,d}$ of
$\overline{\rho}_{x,d}$ for all pairs $(x,d)\in F_\alpha$.
In order to construct $Q_\alpha$ pick a countable and dense subset $Z_\alpha$
of $N_\alpha$. Extend $Z_\alpha$ in countably many steps to close it
under the application of the mappings
$$\overline{\varphi}_{x,d}\text{ and }
\overline{\varphi}^{-1}_{x,d}\text{ for }(x,d)\in E_\alpha$$
and
$$\overline{\rho}_{x,d}\text{ and }
\sigma_{x,d}\text{ for }(x,d)\in F_\alpha$$
while at the same time rendering this set suitable for all t.e.r.s
under consideration. 
Finally choose a bijection
$j:Q_\alpha\to\omega(\alpha+1)\setminus\omega\alpha$ and 
extend $<_T$ in the obvious way to define $T_\alpha$.

\subsection{Extend $\equiv_i$ for successor $i$}

Next we define $\equiv_{i}$ on $T_\alpha$ in the case, that  $i=j+1$.
(Recall that for limit $i=i(\alpha)$ the relation $\equiv_i$ on $T_\alpha$ is
determined by its behaviour below.)
Letting $\alpha^-$ be the maximum of $\alpha\cap C_{i}$, $t\in T_\alpha$ and
$r:=t \uhr \alpha$,
we easily see that we have an infinite set
$\hat{r}\cap t/\!\!\equiv_j$
(here we view $T_\alpha$ as a dense subset of $[T\uhr\alpha]$)
which we partition into $\aleph_0$ infinite sets
$P(t/\!\!\equiv_j,r,n)$, each of them dense in $\hat{r}\cap t/\!\!\equiv_j$.
Note, that the definition of $P$ has to be independent from the choice of $t$.
Finally we define $\equiv_{i}$ on $T_\alpha$ by letting $s\equiv_{i} t$ if and only if
$$s\uhr\alpha^-\equiv_{i} t\uhr\alpha^-\text{ and }(\exists n\in\omega)
s\in P(s/\!\!\equiv_j,s\uhr\alpha^-,n)\text{ and }t\in P(t/\!\!\equiv_j,t\uhr\alpha^-,n).$$
Then by construction, $\equiv_{i}$ is $\infty$-nice as well as $\equiv_{j/i}$.

\subsection{Check isomorphisms}

To finish the induction step,
we verify that the extension procedure for the tree isomosphisms
succeeds.
We only check the most complicated case, that of $\psi_{x,d}$ for
$(x,d)\in F_\alpha$.
We start and show that our choice of $T_\alpha$ admits extensions of the tree
isomorphisms $\psi_{x,d}$ for all $(x,d)\in F_\alpha$.
Fix $(x,d)\in F_\alpha$ and set $c'=(c(x)\cap d)\cup\{\alpha\}$.
First of all extend $\equiv_x$ to $T_\alpha$ by letting $s\equiv_x t$ if and only if 
$s\uhr\gamma\equiv_x t\uhr\gamma$ for all $\gamma\in c(x)$.
(This is the only t.e.r. extending $\equiv_x$ on $T_\alpha$
because of the normality requirement in the definition of t.e.r.)

We can identify $T_\alpha$ with $Q_\alpha$ via the bijection $j$.
This in mind, we show that the unique extension $\psi$
of $\psi_{x,d}$ to $(T\uhr c')/\!\!\equiv_x$ is
an isomorphism onto $T\uhr c'$:
for $s\in T_\alpha$ define
$$\psi'(s/\!\!\equiv_x):=\text{the unique }t\in T_\alpha\text{ with }
\psi_{x,d}(s\uhr\gamma/\!\!\equiv)=t\uhr\gamma\text{ for all }\gamma\in c'\setminus\{\alpha\}$$
and let $\psi=\psi_{x,d}\cup\psi'$.
We check the soundness of this definition.
For $s\in T_\alpha$ there is
$\psi(s/\!\!\equiv_x)=j\circ\overline{\rho}_{x,d}\circ j^{-1}(s)$.
The normality of $T/\!\!\equiv$ guarantees that
the definition of $\psi(s/\!\!\equiv_x)$ is
independent from the choice of the representative $s$.
So it is clear that $\psi$ is indeed a tree homomorphism.

Now, by the definition of $\equiv_x$ on $T_\alpha$,
if $s\not\equiv t$ for $s,t\in T_\alpha$, then there is
some $\gamma\in c'\setminus\{\alpha\}$ with
$s\uhr\gamma\not\equiv_x t\uhr\gamma$ and so
$\psi(s/\!\!\equiv_x)\neq\psi(t/\!\!\equiv_x)$.
On the other hand, for every $s\in T_\alpha$ we have
$$\psi(j\circ\sigma_{x,d}\circ j^{-1}(s)/\!\!\equiv_x)=s.$$
So $\psi$ is bijective also on the top level
of $(T\uhr c')/\!\!\equiv$
and therefore a tree isomorphism.

For reference in later induction steps
we denote this tree isomorphism $\psi$ by
$\psi_{x,d}$.
This completes the inductive construction.

\subsection{Chain homogeneity}

It remains to show that the above construction yields a Souslin tree 
$T$ whose regular open algebra $\B$ is chain homogeneous.
We omit the standard argument proving that $T$ is Souslin.

If we are given two generating
chains $K,K'$ of $\B$ let $X\subset\omega_1$ be a code with respect to $g$
for a pair of wipers on the club
$C\subset\{\alpha:\omega\alpha=\alpha\}$ representing $K$ and $K'$.
Let $D$ be a club set in $\omega_1$ associated to $X$ by $\diaplus$.
Then for each
$\alpha\in E=\{\gamma\in C\cap D\mid\gamma=\bigcup(\gamma\cap C\cap D)\}$
the construction of the tree gives us a tree automorphism
$\varphi_{X\cap\alpha,D\cap\alpha}$ of $T\uhr(C\cap D\cap\alpha)$,
and the union of that increasing chain,
$$\varphi=\bigcup_{\alpha\in E}\varphi_{X\cap\alpha,D\cap\alpha},$$
extends to an automorphism $\varphi$ of $\B$ that carries $K$ to $K'$.

Now let $\A$ be a complete and atomless subalgebra of $\B$
represented by the decent t.e.r. $\equiv$ on $T\uhr C$.
Let $D$ be a club, such that the $\diaplus$-sequence guesses
$D\cap\alpha$ and $X\cap\alpha$ for
a code $X$ for $\equiv$ for all $\alpha\in D$.
Then the construction yields a tree isomorphism
$\psi:(T\uhr C\cap D)/\!\equiv\,\to T\uhr C\cap D$ and therefore a Boolean isomorphism between $\B$ and $\A$.

So let $K$ be a maximal chain of $\B$, $\A=\langle K\rangle^\cm$ and $\psi:\B\cong\A$.
Now $K$ is isomorphic to the generating chain ${\psi^{-1}}''K$
and is thus of that unique order type.

\subsection{Verification of the properties of the $\infty$-nice chain}

Let $i$ be a countable limit ordinal.
Then, by our definition of $\equiv_i$ (in the limit case),
for every node $s\in T\uhr C_i$, its $\equiv_i$-class
is just the intersection over the family $(s/\!\!\equiv_j)_{j<i}$.
This shows that $\A_i=\big\langle \bigcup_{j<i}\A_j\rangle^\cm$,
i.e. continuity of the chain.

For the proof of {\it(b.ii)} let $(\C_\nu\mid\nu<\lambda)$ be given.
By chain homogeneity of $\B$ we can choose $\varphi_0:\A_0\cong\C_0$,
and then inductively extend the given $\varphi_\nu:\A_\nu\to\C_{\nu}$
to $\varphi_{\nu+1}:\A_{\nu+1}\to\C_{{\nu+1}}$
by virtue of condition {\it(a)} as follows.
Choose an isomorphims $\varphi^0_{\nu+1}:\A_{\nu+1}\to\C_{\nu+1}$ and extend
$\psi':=({\varphi^0_{\nu+1}}^{-1}\circ \varphi_\nu)$ by virtue of {\it(a)} to
an automorphism $\psi$ of $\A_{\nu+1}$. Then set
$\varphi_{\nu+1}:=\varphi^0_{\nu+1}\circ\psi$. 

At limit stage $\mu$, the argument given above shows that
$$\left\langle\bigcup_{\nu<\mu}\C_{\nu}\right\rangle^\cm=\C_\mu=
\left\langle\bigcup_{\nu<\mu}{\varphi_\nu}''\A_\nu\right\rangle^\cm.$$
So there is a unique extension $\varphi$ of $\bigcup_{\nu<\mu}\varphi_\nu$
with domain $\A_\mu=\langle\bigcup_{\nu<\mu}\A_nu\rangle^\cm$ and range
$\C_\mu=\langle\bigcup_{\nu<\mu}\C_\nu\rangle^\cm$. 

This finishes the proof of Theorem \ref{thm:hom_for_max}.\hfill$\square$

\section{Features, remarks and open problems}
\label{sec:features}

\subsection{Features of chain homogeneous Souslin algebras}

It is not hard to see that by omission of the additional steps (for properties
a) and b) of Theorem \ref{thm:hom_for_max}) in the
construction we arrive at a small Souslin algebra which is plainly chain
homogeneous. The following two propositions state properties shared by all
chain homogeneous Souslin algebras.

\begin{prp}\label{prp:properties}
A small and chain homogeneous Souslin algebra $\B$
is homo\-geneous in the following strong sense.
For every pair $\A_0,\A_1$ of Souslin sub\-algebras and
$x\in\A_0$ and $y\in\A_1$ where $0<_\B x,y<_\B 1$, there are $2^{\aleph_1}$ distinct
isomorphisms $\varphi:\A_0\to\A_1$ with $\varphi(x)=y$.
\end{prp}
\begin{proof}
We first argue that $\B$ is weakly homogeneous.
$\B$ is assumed to be small.
So for every pair $a,b$ of non-zero elements of $\B$ there are generating
chains $K,K'$ of $\B$ with $a\in K$ and $b\in K'$.
Each isomorphism $\varphi$ between
$K$ and $K'$ satisfies $\varphi(a)\cdot b\neq 0$ and extends to an automorphism of $\B$.

By a theorem proved independently by Koppelberg and by Solovay
(cf. \cite[Thm.18.4.1]{stepanek-rubin}),
every complete and weakly homogeneous Boolean algebra
is a power of a homogeneous factor,
which in our case is isomorphic to $\B$,
because of the c.c.c. satisfied by $\B$.

As there are $2^{\aleph_1}$ distinct generating chains,
chain homogeneity implies that $\B$ has $2^{\aleph_1}$ automorphisms.
The rest is then routine.
\end{proof}

\begin{prp}\label{prp:no-ind-subalg}
\begin{enumerate}[a)]
\item If $\A$ is an atomless, complete subalgebra
of the Souslin algebra $\B$,
and $\A$ and $\B$ are isomorphic,
then no atomless, complete subalgebra of $\B$
can be independent from $\A$.
\item A chain homogeneous Souslin algebra
has no independent pair of
atomless and complete subalgebras.
\end{enumerate}
\end{prp}
\begin{proof}
The proof of  a) is straight forward applying the fact that
a Souslin algebra cannot have anindependent pair of mutually isomorphic
Souslin subalgebras.
For the proof of b) note, that here we have not assumed $\B$ to be small.
In the case of a small Souslin algebra, an application of part a)
suffices to prove b).
But also a big Souslin algebra $\B$ has only maximal chains,
that completely generate subalgebras which are small Souslin algebras.
So assume that $\A_0,\A_1$ form an independent pair of
complete and atomless subalgebras of $\B$.
Then there are maximal chains $K_0\subset\A_0$ and
$K_1\subset\A_1$, which are isomorphic to each other.
The isomorphism between the chains extends to an isomorphism
between the two subalgebras that are completely generated by the chains:
$$\C_0:=\langle K_0\rangle^\cm\subset\A_0\text{ and }\C_1:=\langle
K_1\rangle^\cm.$$ 
But then again, we have an isomorphic pair of subalgebras $\C_0,\C_1$ of $\B$
that cannot be independent unless $\B$ fails to satisfy the countable chain
condition.
\end{proof}

\subsection{The number of non-isomorphic maximal chains}

Classical questions about Souslin trees often circulated around the number
isomorphism types or the number of automorphisms. In this respect, it is worth
noting, that a Souslin algebra which is not chain homogeneous has at least
$2^{\aleph_0}$ non-isomorphic maximal chains and that the regular open
algebras of most Souslin trees considered in set theoretic literature, such as
strongly homogeneous or free Souslin trees have $2^{\aleph_1}$ non-isomorphic
maximal (even generating) chains.

\subsection{On the hypothesis}

Concerning the hypothesis ($\diaplus$) met for our constructions,
we remark the following.
In \cite[\S 6]{abraham-shelah} and \cite[\S 4]{abraham-shelah_aronszajn_trees} a
model of {\sf ZFC + $\neg$SH} is constructed, 
in which there is no homogeneous Souslin tree.
Yet a (chain) homogeneous Souslin algebra
always has a homogeneous Souslinization.
So we have found a {\sf ZFC}-model where Souslin's hypothesis
fails, but which has no chain homogeneous Souslin algebras.

It is open whether we could make do with strictly less than $\diaplus$,
e.g., if the assumption of $\dia$ is sufficient to guarantee
the existence of a chain homogeneous Souslin algebra, cf. the open problems
section below. 

However, we can show that the existence of chain homogeneous Souslin algebras
with properties as strong as stated in Theorem \ref{thm:hom_for_max} cannot be
obtained from $\diamondsuit$ alone.

\begin{prp}
Let $\B^*$ be the big and chain homogeneous algebra of Theorem \ref{thm:big}.
Forcing with $\B^*$ turns any ground model Souslin tree $T$,
which is regularly embeddable into $\B^*$, into a Kurepa tree.
\end{prp}
\begin{proof}
Let $G\subset\B_{\omega_2}$ be a $V$-generic filter, and
let $\A:=\langle T\rangle^\cm$.
For $\alpha<\omega_2$ fix in $V$ an isomorphism
$\varphi_\alpha:\B_\alpha\cong\A$ and in $V[G]$ define
$G_\alpha:=\B_\alpha\cap G$.

Staying in $V[G]$, for $\alpha<\omega_2$, the set
$$b_\alpha:=\{t\in T\mid(\exists g\in G_\alpha\})
\varphi_\alpha(g)\leq_\A t\}$$
is an $\omega_1$-branch of $T$,
and we have $V[G_\alpha]=V[b_\alpha]$ for all $\alpha<\omega_2$.
Finally, if $\gamma<\beta<\omega_2$ then
$b_\beta\in V[b_\beta]\setminus V[b_\gamma]$,
because the map $\varphi_\gamma$ is in $V$.
In particular we have $b_\beta\ne b_\gamma$.
\end{proof}

Since the Iteration Theorem \ref{thm:big}
only assumes the strong properties
of the Souslin algebra of Theorem \ref{thm:hom_for_max}
we have the following result on the hypotheses used.

\begin{cor}\label{cor:hypo}
Assume that there is a Mahlo cardinal.
Then there is a model of {\sf ZFC} + $\diamondsuit$ in which there is no Souslin algebra
with the properties stated in theorem \ref{thm:hom_for_max}.
\end{cor}
\begin{proof}
In \cite{GKH} Jensen considered the generic Kurepa hypothesis
{\sf GKH}:= ``there is a c.c.c. forcing that forces {\sf KH} in the generic extension''.
He shows that if $\kappa$ is a Mahlo cardinal then {\sf GKH} is false in the
Levy-style generic extension collapsing $\kappa$ to become $\aleph_2$.
This partial order always forces $\diamondsuit$ (cf.~\cite[Exercises VIII.J.5/6]{kunen}). 
Finally, $\B_{\omega_2}$ is c.c.c. and forces, by the last Proposition,
{\sf KH}.
\end{proof}

As far as we know, the proofs of Theorems \ref{thm:hom_for_max} and
\ref{thm:big} give the first construction of a big Souslin algebra assuming 
the principle $\diaplus$ only. Jensen's constructions use $\diamondsuit$ and
$\square$, cf. \cite{devlin-johnsbraten} or \cite{GKH}. 
In \cite[\S 5]{jech_mess} gives a forcing that adjoins a so-called Souslin
mess. This is a partial order of partial functions generalising the notion of
a Souslin tree. The regular open algebra of a large enough Souslin mess is a
big Souslin algebra. Laver has constructed a Souslin mess,
only using $\diamondsuit$ and Silver's principle {\sf W}.
({\sf W} is a strengthening of Kurepa's Hypothesis {\sf KH}; cf. \cite[(24.16)]{jech}.) 
% We do not know yet, if and how the principles $\diamondsuit+$ {\sf W} and
% $\diaplus$ are correlated. 

Nevertheless, the following remains open.

\begin{question}
Assuming $\dia$, is there a chain homogeneous Souslin algebra?
\end{question}
This is only a sample question, as $\dia$ could be replaced by any
of its variants known to be strictly weaker than $\diaplus$.

In \cite[\S V.3]{devlin-johnsbraten}, a Souslin tree
with at least $\aleph_2$ automorphisms (a property in common with chain
homogeneous Souslin algebras) is constructed
under $\diaplus$, and on p. 51 the authors remark:
\emph{It is doubtful whether this is provable from $\dia$.}
% It seems very likely,
% that the regular open algebra of that Souslin tree
% is also chain homogeneous as, for instance,
Though it is not hard to see that the regular open algebra of this tree is
homogeneous for generating chains, we do not know whether it is fully chain
homogeneous. 

As argued above, it is impossible to carry out our construction of
Section \ref{sec:mainproof} under $\dia$ alone. 
% Nevertheless, we ask:

% This is reinforced by the fact that
% the small Souslin algebra constructed in Section \ref{sec:big}
% has stronger homogeneity properties
% of which we can prove, that they do not exist under
% the assumption of $\dia$ alone (cf. Corollary \ref{cor:hypo}).
% And the methods used to implement these stronger homogeneity properties strongly
% resemble those used in the proof
% of Theorem \ref{thm:hom_for_max}.

% So we join in the doubts of Devlin and Johnsbr\aa ten
% cited above and conjecture
% that there is a model of {\sf ZFC} + $\dia$ in which
% there are no chain homogeneous Souslin algebras.

\subsection{Souslin's hypothesis minus one}

\begin{question}
Is there a model of {\sf ZFC} with exactly one Souslin
line (up to isomorphism)?
\end{question}
Such a property could be paraphrased as ``Souslin's hypothesis minus one'' as
there would be only one counter-example, essentially, i.e. only one without
separable intervals.
In such a model the associated Boolean algebra would have to be chain
homogeneous. Anyway, also if we refocus the question on Souslin trees or
Souslin algebras, then the respective algebra must admit isomorphisms to all
of its Souslin subalgebras, a feature in common with the chain homogeneous
Souslin algebras. 

In an attempt to construct such a model by variations of the well-known
iterated forcing techniques of Solovay-Tennenbaum (variation for a strongly
homogeneous Souslin tree by Larson,
cf. {\cite[\S 4]{larson}}) or that of Jensen (variation for a free Souslin tree by
Abraham, cf. {\cite[\S4]{abraham-shelah_aronszajn_trees}}), 
which leave one a priori specified Souslin tree intact, the main
obstacle is to ensure that enough of the homogeneity properties of the
preserved Souslin algebra remain valid in the final model.

\subsection*{Acknowledgements}
  Thanks are due to Sabine Koppelberg for asking the initial question,
  to Stevo Todor\c{c}evi\'{c} for asking a question which led to a
  significant improvement of my results and to Peter
  Krautzberger and Stefan Geschke for support, remarks and discussions.

\bibliographystyle{alpha}
\bibliography{souslin,buecher}
%\addcontentsline{toc}{chapter}{References}

\end{document}